\documentclass[12pt,fullpage,leqno]{article}

\usepackage{amsfonts}
\usepackage{amsmath}
\usepackage{epsfig}
\usepackage{graphicx}
\usepackage{xcolor} 
\input epsf

\def\m{\mu}

\def\d{\delta}

\def \diam{\hbox{\rm diam }}
\newcommand{\ba}{\begin{eqnarray}}
\newcommand{\ea}{\end{eqnarray}}

\newcommand{\qeda}{\hfill $\square$}

\newcommand{\e}{\varepsilon}

\newenvironment{Proofc}[1]{\smallskip\par\noindent\textsc{#1}\quad}%
  {\hfill$\Box$\bigskip\par}

\nonstopmode

\numberwithin{equation}{section}
\newcommand{\<}{\langle}
\renewcommand{\>}{\rangle}
\newcommand{\tr}{\operatorname{tr}}
\newcommand{\sgn}{\operatorname{sgn}}

\newcommand{\R}{\mathbb{R}}

\renewcommand{\P}{\mathcal{P}}

\newcommand{\Z}{\mathbb{Z}}

\newtheorem{de}{Definition}[section]
\newtheorem{lem}[de]{Lemma}
\newtheorem{prop}[de]{Proposition}
\newtheorem{thm}[de]{Theorem}
\newtheorem{cor}[de]{Corollary}
\newtheorem{rem}[de]{Remark}

\newcommand{\dd}{\partial}
\renewcommand{\d}{d}

\newcommand{\supp}{\operatorname{supp}}

\newcommand{\lap}{\Delta}

\newcommand{\gammab}{{\overline{\gamma}}}
\newcommand{\ms}{M_\textup{sing}}
\newcommand{\ns}{N_\textup{sing}}
\newcommand{\ps}{\mathcal{P}_\textup{sing}}
\newcommand{\n}{\bar N}

\begin{document}

\title{\bf Pointwise regularity of the free boundary\\ 
for the parabolic obstacle problem}
 \author{Erik Lindgren\footnote{{erik.lindgren@math.ntnu.no, Dept. of Mathematical
Sciences, NTNU, 7491 Trondheim, Norway}} and 
R\'egis Monneau\footnote{Universit\'e Paris-Est, CERMICS, Ecole des Ponts ParisTech,
   6 et 8 avenue Blaise Pascal, Cit\'e Descartes, Champs-sur-Marne,
   77455 Marne-la-Vall\'ee Cedex 2}}

\maketitle

\vspace{20pt}



 \centerline{\small{\bf{Abstract}}}
\noindent {\small{We study the parabolic obstacle problem
$$\lap u-u_t=f\chi_{\{u>0\}}, \quad u\geq 0,\quad f\in L^p \quad \mbox{with}\quad f(0)=1$$
and obtain two monotonicity formulae, one that applies for general free boundary points and one for singular free boundary points. 
These are used to prove a second order Taylor expansion at singular points (under a pointwise Dini condition),
with an estimate of the error (under a pointwise double Dini condition).
Moreover, under the  assumption that $f$ is Dini continuous,
we prove that the set of regular points is locally a (parabolic) $C^1$-surface and that the set of singular points is locally contained in a union of (parabolic) $C^1$ manifolds.}}\hfill\break

 \noindent{\small{\bf{AMS Classification:}}} {\small{35R35.}}\hfill\break
 \noindent{\small{\bf{Keywords:}}} {\small{Obstacle problem, free boundary, parabolic equation, monotonicity formula, singular set, Dini condition.}}\hfill\break


\tableofcontents

\section{Introduction}

We study the parabolic obstacle problem
\begin{equation}\label{eq:obstacle}
\left\{\begin{array}{l}
\left.\begin{array}{l}\lap u-u_t=f\chi_{\{u>0\}} \\
\\
u\geq 0
\end{array}\right|
\textup{ in $\R^n\times (-1,0]$},\\
\\
0\in \left\{u=0\right\}\cap \partial\{u>0\},\\
\\
u,f\in L^p(\R^n\times (-1,0])\textup{ and }f(0)=1,\\
\\
\supp u\subset Q_1^-:= B_1\times (-1,0].
\end{array}\right.
\end{equation}
Here $B_1=\left\{x\in\R^n,\ |x|<1\right\}$ is the unit ball.
By $\chi_{\{u>0\}}$ we denote the
characteristic function of the set $\{u>0\}$ and $\left\{u=0\right\}\cap \partial\{u>0\}$ is referred to as the free boundary.
Notice that we do not simply define the free boundary as $\partial \left\{u>0\right\}$, in order to 
allow $0\in Q_1^-$ to be a point of the free boundary in the case where $u>0$ on $Q_1^-\backslash \left\{0\right\}$ with $u(0)=0$,
and to exclude $\left(\overline{B}_1\times \left\{0\right\}\right)\backslash \left\{0\right\}$ from the free boundary.
Moreover, we consider
exponents $p\in (1,\infty)$ and we assume that the origin is a
Lebesgue point of $f$ in order for $f(0)$ to be well defined, whenever
necessary.

We note that for any local solution  (i.e., a solution in $B_1\times (-1,0]$) 
we can obtain a solution of \eqref{eq:obstacle} with slightly different $f$ by
multiplying with a suitable cut-off function.

The parabolic obstacle problem arises naturally in many different contexts
such as in the modelling of ice-melting (the Stefan problem) 
and in pricing of American options. But it has also an interest of its own.

The present paper can be seen as a continuation of the study commenced in \cite{LM11}, 
where mainly the regular points (see the next section) of the free boundary was studied. 
We develop similar techniques for singular points and 
perform a finer analysis in the case when the function $f$ 
in the right hand side of \eqref{eq:obstacle} is assumed to be Dini continuous.

\subsection{Main results} 

In order to present the results of the paper, we define two different moduli of continuity. 
First, the usual
$L^p$ modulus of continuity (which was used in \cite{LM11}) with $Q_\rho^-=B_\rho\times (-\rho^2,0]$:
$$\tilde\sigma_p^f(r)=\sup_{\rho\in(0,r]}\left(\frac{1}{|Q_\rho^-|}\int_{Q_\rho^-}
  |f(x,t)-f(0)|^p\d x\d t\right)^\frac{1}{p},$$
which is finite for instance if $0$ is a Lebesgue point for $f$.
Second, the modulus of continuity in $L^p_\gamma$-average
$$\sigma_p^f(r)=\sup_{\rho\in(0,r]}\left(\frac{1}{\rho^2}\int_{S_\rho}
  |f(x,t)-f(0)|^p\d \gamma  \right)^\frac{1}{p},$$
where $S_\rho=\R^n\times (-\rho^2,0]$ and  $\d \gamma = G(x,-t)\d x\d t$ 
where $G(x,t)$ is the Heat kernel. Whenever there is no possible confusion, $\sigma_p^f$ will be simply written as
$\sigma_p$, $\sigma^f$ or even $\sigma$, and similarly $\tilde{\sigma}_p^f$ will be simply written as
$\tilde{\sigma}_p$, $\tilde{\sigma}^f$ or even $\tilde{\sigma}$. 
Furthermore, a modulus of continuity $\sigma(r)$ is said to be {\it Dini} if
\begin{equation}\label{eq:dini}
\int_0^1\frac{\sigma(r)}{r}\d r<\infty,
\end{equation}
and a function is said to be Dini continuous if its (usual) modulus of continuity satisfies \eqref{eq:dini}. 
We say that $\sigma$ is {\it double Dini} if
\begin{equation*}
\int_0^1  \frac{1}{r} \left(\int_0^r\frac{\sigma(s)}{s}\ \d s\right)\d r<\infty.
\end{equation*}

\begin{rem}{\bf (Relation between $\tilde{\sigma}$ and $\sigma$)}\\ 
As can be seen in Appendix A, taking this quite unusual
  definition of the modulus of continuity gives actually a weaker assumption
  than assuming anything for the usual $L^p$ modulus of continuity. In
  fact, from Proposition \ref{prop:modcont}, we have the following relation for $p\ge 1$:
$$\sigma(r)\leq C\left(\tilde{\sigma}(r)+\tilde{\sigma}(\sqrt{ r})+\tilde{\sigma}(1) e^{-\frac{c}{r}}\right).$$
In particular, $\sigma$ is Dini (resp. double Dini) if $\tilde{\sigma}$ is Dini (resp. double Dini).

One can easily give examples of functions which are Dini in the sense of $\sigma$ but not in the sense of $\tilde \sigma$. For instance, consider a function of the form
$$
f(x,t)=g(t)h\left(\frac{x}{-t}\right), 
$$
where $g$ is $C^\alpha$ with $\alpha\in (0,1)$, $g(0)=0$ and 
$$
\int_{\R^n} h(y) e^{-\frac{y^2}{4}} dy <\infty.
$$
Then $\sigma_1^f$ will be Dini while for appropriate choices of $g$ and $h$, $\tilde \sigma_1^f$ will not even be finite. Take for instance take $g=|t|^\alpha$ and $h(y)=y^4$.
\end{rem}

The first result is that we are able to classify the possible types of free boundary points as seen
below. We introduce the notion of \emph{blow-ups}.
Let $u$ be a solution of \eqref{eq:obstacle} and $X_0=(x_0,t_0)$ be a free boundary point. 
If we can find a subsequence $u_{X_0,r_j}$ of
$$u_{X_0,r}(x,t)=\frac{u(rx+x_0,r^2t+t_0)}{r^2},$$
for $r_j\to 0$, converging locally uniformly to a limit $u_0$, 
then we say that $u_0$ is a \emph{blow-up} of $u$ at $X_0$.

\begin{thm}\label{thm:class}{\bf (Classification of free boundary points via the Weiss functional)}\\
Let $p> (n+2)/2$ with $p\ge 2$ and let $u$ be a solution of \eqref{eq:obstacle}. Suppose further that $\tilde{\sigma}_p(1)$ is finite and 
that $\sigma_p$ is {\bf Dini}. 
Then the function
$$E(r,u)=\frac{1}{r^4}\int_{S_r}\left(|\nabla u|^2+2u+\frac{u^2}{t}\right)\d\gamma$$
has a limit at $r=0$, denoted by $E(0^+,u)$. Moreover, one of the following alternatives holds
\begin{enumerate}
\item $E(0^+,u)=0$ and the origin is said to be a \emph{degenerate point},
\item $E(0^+,u)=15/2$ and the origin is said to be a \emph{regular point},
\item $E(0^+,u)=15$ and origin is said to be a \emph{singular point}.
\end{enumerate}
\end{thm}
We define
$$
N_{\textup{sing}}(u,\rho)=\inf_{P\in \P_{\textup{sing}}}
\left(\frac{1}{\rho^4}\int_{S_\rho} |u-P|^2\d\gammab\right)^\frac{1}{2}$$
and 
$$ M_{\textup{sing}}(u,\rho)=\sup_{r\in (0,\rho]} N_{\textup{sing}}(u,\rho),$$
where $\displaystyle \d\gammab = \frac{G(x,-t)}{-t}dxdt$ and
\begin{equation}\label{eq::el1}
\P_{\textup{sing}}=\left\{\begin{array}{lr}
P(x,t)=\frac12 {}^tx\cdot Q\cdot x +mt,\\P\geq 0 \quad \mbox{on}\quad \R^n\times \R^-,\\\lap P -P_t =1,\\m\in [-1,0].\end{array}\right\}.
\end{equation}
As can be seen in Lemma 6.3 in  \cite{CPS04}, $\ps$ is the set of possible blow-up solutions at singular points.

At singular points we obtain results similar to those at regular points (cf. \cite{LM11}),
i.e., an explicitly controlled Taylor expansion of second order.

\begin{thm}\label{thm:mainsing}{\bf (Modulus of continuity at singular free boundary points)}\\
Let $p>(n+2)/2$ with $p\ge 2$. Assume that $u$ satisfies (\ref{eq:obstacle}), that $\sigma_2$ is {\bf double Dini} and that $\tilde \sigma_p(1)<\infty$.
Then there exist $\alpha\in (0,1]$ and constants
  $C>0,M_0,r_0\in (0,1)$ such that 
$$M_{\textup{sing}}(u,r_0)\leq M_0$$
implies the existence of $P_0\in \P_\textup{sing}$ satisfying for all $r\in (0,r_0)$
$$\left(\frac{1}{r^4}\int_{S_r} |u-P_0|^2\d\gammab\right)^\frac{1}{2}\leq
C\left(M_{\textup{sing}}(u,r_0)r^\alpha+\int_0^r\frac{\Sigma_p(s)}{s}\d
  s+r^\alpha\int_r^1\frac{\Sigma_p(s)}{s^{1+\alpha}}\d s\right),$$
where
$$\Sigma_p(\tau)=\sigma_p(\tau)+\int_0^\tau
  \frac{\sigma_p(r)}{r}.$$
\end{thm}

What is not so usual in the literature, except in \cite{Mon09} and \cite{LM11} is that the result above is completely pointwise.

Before presenting the rest of our results, let us introduce some notation related to the parabolic distance.

\begin{de}\label{def:c1} {\bf ((parabolic) $C^1$ functions)}\\
Let $d(x,t)=\sqrt{|x|^2+|t|}$ denote the parabolic distance. 
We say that $f:E\subset\R^{n+1}\to \R$ is $C^1$ with
respect to $d$, with derivative $(g,0)$ if there exists a modulus of continuity $\omega$ (i.e. with $\omega(0^+)=0$)
such that for all $X, X+H \in E$ with $H=(h^x,h^t)$ there holds
$$\left\{\begin{array}{l}
|f(X+H)-f(X)-g(X)h^x|\le d(H)\omega(d(H)),\\
|g(X+H)-g(X)|\le \omega(d(H))
\end{array}\right.$$
When there is no possible confusion, we will simply
say $C^1$ when we mean $C^1$ with respect to $d$.
\end{de}

\begin{rem} 
Any $C^1$-function $f:\R^{n+1}\to \R$ (in the euclidean setting) is $C^1$ with respect to $d$, 
with derivative $(\nabla_{\R^n}f,0)$.
\end{rem}

\begin{de}\label{defi::C1manif} {\bf ((parabolic) $C^1$-manifold)}\\
The notion of $C^1$ in Definition \ref{def:c1} induces naturally a notion of 
(parabolic) $(k+1)$-dimensional $C^1$-manifolds in $\R^{n+1}$ as graphs of
maps from $\R^{k+1}\to \R^{n-k}$, where each coordinate is a
$C^1$ function as in Definition \ref{def:c1}.
\end{de}

Given a point $X_0$ of the free boundary $\left\{u=0\right\}\cap \partial \left\{u>0\right\}$, 
if $f(X_0)>0$, we can consider the values of $E\left(0^+, \frac{u(X_0+\cdot)}{f(X_0)}\right)$.
In view of Theorem \ref{thm:class}, when $\sigma^{f(X_0+\cdot)}_p$ is Dini, we can decide if $X_0$ is degenerate, regular or singular.
We can then split the free boundary into the three associated sets:
$$\left\{u=0\right\}\cap \partial \left\{u>0\right\}=\Gamma_d\cup \Gamma_r\cup \Gamma_s$$
where $\Gamma_d$ is the set of degenerate points, $\Gamma_r$ the set of regular points and $\Gamma_s$ the set of singular points.

We say that $f$ is {\it uniformly Dini continuous} in an open set ${\mathcal Q}\subset \R^{n+1}$ if there exists a Dini modulus of continuity $\bar \sigma$ such that
$$|f(X)-f(Y)|\le \bar \sigma(d(X-Y)) \quad \mbox{for all}\quad X,Y\in {\mathcal Q}.$$
For such $f$, when we assume moreover that $f>0$ on ${\mathcal Q}$, the degenerate points are impossible (see Proposition \ref{prop:nondeg2}).
Therefore the free boundary only splits in two parts 
$\Gamma_r$ and $\Gamma_s$, which consists of the regular and singular points, respectively.

In the case of regular points, we have the result below, which is a by-product of Theorem \ref{thm:class} 
in the present paper and Theorem 1.7 in \cite{LM11}.

\begin{thm}{\bf (Regularity of the regular part of the free boundary)}\label{thm:mainreg}\\
Consider a solution $u$  of (\ref{eq:obstacle}), and assume that $f$ is  {\bf uniformly Dini continuous} in a neighborhood of the origin. 
Then in a neighborhood of the origin, the set $\Gamma_r$  of regular free boundary points is an open subset of the free boundary $\Gamma := \left(\partial\left\{u>0\right\}\right)\cap \left(\partial \left\{u=0\right\}\right)$ and 
around any regular free boundary point there exists a neighborhood $V$ such that 
$V\cap \Gamma$ is locally a $C^1$ hypersurface with respect to the parabolic distance. 
More precisely, up to a rotation of the spatial coordinates
$$V\cap \Gamma = \left\{(x,t) \quad \mbox{such that}\quad x_n=\tilde{f}(x',t) \quad \mbox{with}\quad  (x',t)\in V'\right\},$$
where $x'=(x_1,...,x_{n-1})$, the set $V'$ is an open set in $\R^n$, and $\tilde{f}: V'\to \R$ is a $C^1$ function.
\end{thm}

If $f$ is more regular it is reasonable to expect higher regularity of the free boundary at regular points. For instance, if $f$ is Lipschitz in both $x$ and $t$, it seems to be possible to adapt the proof of Proposition 11.1 in \cite{SUW09} to our situation in order to obtain that the free boundary is also Lipschitz with respect to $t$.

If the origin is a singular point of the free boundary, with $\tilde{\sigma}_p(1)$ finite and $\sigma_p$ Dini for some $p>(n+2)/2$ with $p\ge 2$, 
it is possible to show the uniqueness of the blow-up limit (see Corollary \ref{cor:unique}). 
Hence we can introduce the following notation, 
which is needed to state the next main result of the paper.

\begin{de}\label{defi::el15}{\bf (Singular set)}\\
Given a singular point of the free boundary, we consider its unique blow-up limit $P(x,t)=\frac12 {}^tx\cdot Q \cdot x + mt$ with $Q\ge 0$, $m\le 0$ and $\Delta P-P_t>0$. Morever, we denote by $\Gamma(k+1)$ the set of singular free boundary points such that
the matrix $Q$ has a kernel of dimension $k\in \left\{0,...,n-1\right\}$ and by $\Gamma(n+0)$ the set of singular points such that $m<0$.
\end{de}
Notice that by Definition \ref{defi::el15}, we have
$$\Gamma_s = \Gamma(n+0)\cup \left(\bigcup_{k=0,...,n-1}\Gamma(k+1)\right).$$

\begin{thm}\label{thm:mainsing2}{\bf (Structure of the singular set of the free boundary)}\\
Let $u$ satisfy (\ref{eq:obstacle}) and suppose that $f$ is {\bf uniformly Dini continuous}
in a neighborhood of the origin. Then
\begin{enumerate}
\item For each $k\in \left\{0,...,n-1\right\}$, the set $\Gamma(k+1)$ is locally contained in a $(k+1)$-dimensional (parabolic) $C^1$-manifold in $\R^{n+1}$.
\item $\Gamma(n+0)$ is locally contained in a (euclidean) $C^2$-graph of the form $t=t(x_1,\ldots,x_n)$.
\item $\Gamma_s$ is locally contained either in an (euclidean) $C^2$-graph of the form $t=t(x_1,\ldots,x_n)$
or in a $((n-1)+1)$ dimensional $C^1$-manifold.
\end{enumerate}
\end{thm}

\begin{rem} 
To the authors' knowledge, Theorem \ref{thm:mainreg} is a new result unless we assume $f=1$ 
and part 1 and 3 in Theorem \ref{thm:mainsing2} are new results even in the case when $f$ is constant.
\end{rem}

\subsection{Known results} 

The parabolic obstacle problem has been in the focus of attention for many years. 
The regularity of the solution and also the properties of the free boundary near regular free boundary points, 
are by now fairly well understood. The first paper treating this satisfactory is \cite{CPS04}, 
where the case $f=1$ is treated. Under slightly different assumptions the problem 
has also been studied in for instance \cite{BDM05}, \cite{BDM06}, \cite{EL12} and \cite{LM11}.

When it comes to singular points, the only result the authors are aware of are the ones in \cite{Bla06}, 
where the author proves that a general singular free boundary point lies in a $C^\frac12$-manifold 
(in both time and space), and that the whole set of singular points lies in a union of $C^\frac12$-manifolds (here in the usual euclidean setting). 

This is in contrast with the elliptic case, where also the set of 
singular free boundary points is very well understood. 
The study of the singular part of the free boundary was initiated in \cite{Sch77} and \cite{CR77},
and later in \cite{Caff} it was proved that a general singular free boundary point 
lies in a $C^1$-manifold, and that the whole set of singular points lies in a union of $C^1$-manifolds. 
This was also generalized into a more general setting in \cite{Mon09}. 

In \cite{Wei99}, a monotonicity formula was developed which has showed to be a very strong tool 
for studying free boundary problems in general. In \cite{Mon03}, 
this technique was even further refined where the author introduced a monotonicity formula 
for singular points, which is then used to prove sharp results about the structure of the singular set. 
These two formulas are used in \cite{Mon09} to obtain a pointwise second order 
Taylor expansion around free boundary points. 
Together with \cite{LM11}, the present paper can be seen as the parabolic counterpart to \cite{Mon09}.

\subsection{Organization of the paper} 

The paper is organized as follows. In Section \ref{sec:reg}, we recall certain well known results for parabolic equations and also some regularity estimates obtained in \cite{LM11}. The next section, Section \ref{sec:mon}, is concerned with two monotonicity formulae, one for general free boundary points and one for singular free boundary points. The former is used to classify free boundary points into three different types and the latter is used to prove the uniqueness of blow-ups at singular free boundary points. 
These results are then applied in Section \ref{sec:dini}, where we under the assumption that $f$ in the right hand side of \eqref{eq:obstacle} is Dini, prove certain structural results (Theorems \ref{thm:mainreg} and \ref{thm:mainsing2}) 
about both the set of regular free boundary points and the set of singular free boundary points. 
This is followed by Section \ref{sec:decay} (which is independent of Section \ref{sec:dini}), where we obtain pointwise 
decay estimates (Theorem \ref{thm:mainsing}) for a certain $L^p$-average of $u$ at singular points. 
Finally, the last section is divided in three appendices. In Appendix A, we give an explanation of the relation between the two different moduli of continuity defined in the beginning of the paper, in Appendix B we prove an extension theorem 
\`a la Whitney in the parabolic setting, and in Appendix C we give a short motivation why the integration by parts performed at several stages is well justified.

\subsection{Notation}

Throughout the paper we will use the following notation:
$$\begin{array}{ll}
u_t=\dd_t u =\frac{\partial u}{\partial t}&\textup{- the time derivative}\\
\lap u=\sum_{i=1}^n \frac{\dd^2 u}{\dd {x_i}^2} &\textup{- the Laplace operator}\\
Hu= \Delta u -u_t &\textup{- the heat operator}\\
Q_r(x_0,t_0)=B_r(x_0)\times (t_0-r^2,t_0+r^2) &\textup{- a parabolic cylinder}\\
Q_r^-(x_0,t_0)=B_r(x_0)\times (t_0-r^2,0]&\textup{- a half cylinder}\\
Q_r = Q_r(0,0), \quad Q_r^-=Q_r^-(0,0)
&\textup{- simplified notation}\\
S_r=\R^n\times (-r^2,0]&\textup{- an infinite strip}\\
T_r=\R^n\times (-4r^2,-r^2]&\textup{- another infinite strip}\\
d(x,t)=\sqrt{|x|^2+|t|}&\textup{-  the parabolic distance}\\
\displaystyle G(x,-t)=\frac{c_n}{(-t)^{\frac{n}2}}e^{\frac{x^2}{4t}}&\textup{- the backward heat kernel}\\
Lu=x\cdot \nabla
  u+2tu_t-2u&\textup{- homogeneity operator}\\
\overline{G}(x,-t)=\frac{G(x,-t)}{-t}\\
\d\gamma= G(x,-t)\d x \d t\\
\d\overline{\gamma}= \overline{G}(x,-t)\d x\d t\\
\d\gamma^{s}=G(x,-s)\d x
\end{array}$$
$L^p_\mu$, the weighted $L^p$-space equipped with the norm 
$$\left(\int |f|^p\d \mu\right)^\frac{1}{p}$$
$L^p_t(L^q_x(A\times B))$, the space equipped with the norm 
$$\left(\int_A\left(\int_B f^q\d x\right)^\frac{p}{q}\d t\right)^\frac{1}{p}$$
$W^{2,1}_p(Q_\rho^-)$, the space
endowed with the norm
$$||u||_{W^{2,1}_p(Q_\rho^-)}=||u||_{L^p(Q_\rho^-)}+||\nabla
u||_{L^p(Q_\rho^-)}+||D^ 2u||_{L^p(Q_\rho^-)}+||u_t||_{L^p(Q_\rho^-)}$$
$L^p(I;X)$, the space endowed with the norm
$$\|u\|_{L^p(I;X)}=\| \|u(t)\|_{X}\|_{L^p(I)}$$

\section{Regularity and growth estimates}\label{sec:reg}

Here we present some known results that will be needed in the paper.

\begin{thm}\label{thm:interior}{\bf (Parabolic interior $L^p$-estimates)}\\
Let $p\in (1,\infty)$. If $u\in L^p(Q_r^-)$ and $Hu\in L^p(Q_r^-)$ then there exists $C=C(r,n,p)$ such that
$$||u||_{W^{2,1}_p(Q_{r/2}^-)}\leq C\left(||u||_{L^p(Q_r^-)}+||Hu||_{L^p(Q_r^-)}\right).$$
\end{thm}

The result above is a special case of Theorem 7.22 on page 175 in \cite{Lie96}.

\begin{thm}\label{thm:globalpar}{\bf (Parabolic $L^p$-estimates)}\\
Let $p\in (1,\infty)$. If $u\in L^p(Q_1^-)$ with $\mbox{supp}\ u\subset Q_1^-$ and $Hu\in L^p(Q_1^-)$, then there exists $C=C(n,p)$ such that
$$||u||_{W^{2,1}_p(Q_{1}^-)}\leq C||Hu||_{L^p(Q_1^-)}.$$
\end{thm}

The above result is given in Proposition 7.18 on page 173 in \cite{Lie96}.

\begin{thm}\label{thm:sobolev} {\bf (Parabolic Sobolev embedding)}\\ 
Let $u\in W_p^{2,1}(Q_r^-)$ with $p\in ((n+2)/2,\infty)$. Then there exists $C_*=C_*(r,n,p)$ such that
$$||u||_{C^{\alpha}(Q_r^-)}\leq C_*||u||_{W^{2,1}_p(Q_r^-)},$$
with $\alpha=2-\frac{n+2}{p}$ and where $C^\alpha(Q_r^-)$ refers to the parabolic H\"older space.
\end{thm}

This result is contained in Lemma 3.3 on page 80 in \cite{LSU67}.

\begin{prop}\label{prop:quadgrowth} {\bf (Quadratic growth)}\\
Let $u$ be a solution of \eqref{eq:obstacle} for $p\in ((n+2)/2,\infty)$ with $\tilde{\sigma}_p(1)$ finite.
Then there is $C=C(||f||_{L^p(Q_1^-)},||u||_{L^p(Q_1^-)},\tilde{\sigma}_p(1),n,p)$ such that 
\begin{equation}\label{eq::el6}
\sup_{Q_r^-}|u|+\left(\frac{1}{|Q_r^-|}\int_{Q_r^-}|u|^p\right)^\frac1p \leq Cr^2
\end{equation}
whenever $r\le 1$.
\end{prop}

\noindent {\bf Proof of Proposition \ref{prop:quadgrowth}}\\
The result is contained in Proposition 4.1 and Corollary 4.2 in \cite{LM11} for $r<r_0$ for some $0<r_0<1$.
For $r\in [r_0,1]$, the result follows from the $L^p$ estimate (Theorem \ref{thm:globalpar}) and the Sobolev embedding (Theorem \ref{thm:sobolev}).
\qeda

\begin{cor}\label{cor:lpprim}{\bf (Estimates in $L^{p'}_\gamma$-spaces)}\\ 
Assume the hypotheses of Proposition \ref{prop:quadgrowth} with moreover $p\ge 2$.
Then with $Lu=x\cdot \nabla u+2tu_t-2u$ there is $C=C(||f||_{L^p(Q_1^-)},||u||_{L^p(Q_1^-)},\tilde{\sigma}_p(1),n,p)$ such that
$$\int_{B_1\times (-r^2,0]}\left(|Lu|^{p'}+|u|^{p'}\right)\d \gamma\leq Cr^{2+2p'}$$
and
$$\int_{B_1\times (-r^2,0]}\left\{|\nabla u|^2+|u|+\frac{u^2+(Lu)^2}{-t}\right\}\d \gamma\leq Cr^{4},$$
for $r\le 1$.
\end{cor}

\noindent {\bf Proof of Corollary \ref{cor:lpprim}}\\
\noindent {\bf Step 1: Preliminaries}\\
Define, as earlier, for $0<r\le 1$:
$$v_r(x,t)=\frac{u(rx,r^2t)}{r^2}.$$
Then, from Proposition \ref{prop:quadgrowth}
$$\int_{Q_1^-}|v_r|^p\d x\d t\leq C,$$
and
$$\int_{Q_1^-}|Hv_r|^p\d x\d t\leq C.$$
Therefore, by Theorem \ref{thm:interior}
$$\|v_r\|_{W^{2,1}_p}(Q_\frac12^-)\leq C.$$
In particular, since $p'\leq p$, 
$$\|\nabla v_r\|_{L^{p'}(Q_\frac12^-)}+\|(v_r)_t\|_{L^{p'}(Q_\frac12^-)}\leq C.$$
Scaling back, this implies with $\rho=r/2$:
\begin{equation}\label{eq::el2}
\int_{Q_\rho^-}\left(|\nabla u|^{p'}+\rho^{p'}|u_t|^{p'}\right)\d x\d t\leq C\rho^{p'}\rho^{n+2}
\end{equation}
and
\begin{equation}\label{eq:rprimest} 
\int_{Q_\rho^-}|u|^{p'}\d x\d t\leq C\rho^{2p'}\rho^{n+2}.
\end{equation}
Therefore (\ref{eq::el2}) and (\ref{eq:rprimest}) are true for any $\rho<1/2$.
Notice that (\ref{eq::el2}) and (\ref{eq:rprimest}) still hold for $\rho\in [\frac12, 1]$, because of Theorem \ref{thm:globalpar}.
In order to prove  the desired estimates on $B_1\times (-r^2,0]$, 
we split the estimates respectively on $Q_r^-$ and on $(B_1\backslash B_r)\times (-r^2,0]$.\\
\noindent {\bf Step 2: Estimates on $Q_r^-$}\\
\noindent {\bf Step 2.1: $0\le G(x,-t)\le \frac{c}{(2^{-k}r)^n}$ for $(x,t)\in C_k$}\\
Let $C_k=Q^-_{2^{-k}r}\setminus Q^-_{2^{-k-1}r}$. 
Then for $(x,t)\in C_k$ we have either $t\leq -r^24^{-k-1}$ which implies
$$|G(x,-t)|\leq \frac{C}{(2^{-k} r)^n},$$
or $|x|^2\geq 4^{-k-1}r^2$ which implies (where $C$ is a generic constant that may change from inequality to inequality)
$$|G(x,-t)|\leq C \frac{e^{-\frac{4^{-k-2}r^2}{-t}}}{(-t)^\frac{n}{2}}\leq \frac{C}{(2^{-k} r)^n}$$
since maximum will be attained when $t=-C4^{-k}r^2$.\\
\noindent {\bf Step 2.2: The weighted integrals on $Q^-_{r}$}\\
Using \eqref{eq:rprimest} this gives with $r_k=2^{-k}r$:
\begin{align*}
\int_{Q_{r}^-}\left(|Lu|^{p'}+|u|^{p'}\right)\d\gamma&\leq 
C\int_{Q_{r}^-}\left(|\nabla u|^{p'}|x|^{p'}+|u_t|^{p'}|t|^{p'}+|u|^{p'}\right)\d\gamma\\
&\leq C\sum_{k=0}^\infty \int_{C_k}\left(|\nabla u|^{p'}r_k^{p'}+|u_t|^{p'}r_k^{2p'}+|u|^{p'}\right)\d\gamma\\
&\leq C \sum_{k=0}^\infty (2^{-k}r)^{n+2+2p'}\frac{1}{(2^{-k}r)^n}\\
&\leq C \sum_{k=0}^\infty 4^{-k(1+p')}{r^{2+2p'}}\leq Cr^{2+2p'}.
\end{align*}
By the same reasoning, 
$$\int_{Q_{r}^-}\left\{|\nabla u|^2 +|u|+\frac{u^2+(Lu)^2}{-t}\right\}\d\gamma\leq Cr^4.$$
\noindent {\bf Step 3:  Estimates on $(B_1\setminus B_r)\times (-r^2,0]$}\\
\noindent {\bf Step 3.1:  $0\le G(x,-t)\le C\frac{e^{-c4^{k}}}{r^n}$ for $(x,t)\in \tilde{C}_k$}\\
Similarly we define $\tilde{C}_k=(B_{2^{k+1}r}\setminus B_{2^kr})\times (-r^2,0]$. 
We observe that for $(x,t)\in \tilde{C}_k$ we have for some $c>0$,
$$|G(x,-t)|\leq C\frac{e^{-c4^k}}{r^n},$$
again since the max is attained at $t=-C4^{k}r^2$.\\
\noindent {\bf Step 3.2: The weighted integrals on $(B_1\setminus B_r)\times (-r^2,0]$}\\
As before we obtain by invoking \eqref{eq:rprimest} with $r_k=2^{k+1}r$:
\begin{align*}
\int_{(B_1\setminus B_r)\times (-r^2,0]}\left(|Lu|^{p'}+|u|^{p'}\right)\d\gamma&\leq 
C\int_{(B_1\setminus B_r)\times (-r^2,0]}\left(|\nabla u|^{p'}|x|^{p'}+|u_t|^{p'}|t|^{p'}+|u|^{p'}\right)\d\gamma\\
&\leq C\sum_{2r\le 2^{k+1}r\le 2}
\int_{\tilde{C}_k}\left(|\nabla u|^{p'}r_k^{p'}+|u_t|^{p'}r_k^{2p'}+|u|^{p'}\right)\d\gamma\\&\leq
C \sum_{k=0}^\infty (2^{k+1}r)^{n+2+2p'}\frac{e^{-c4^k}}{r^n}\leq
Cr^{2+2p'}.
\end{align*}
Repeating this for $|\nabla u|^2$, $|u|$,  $\frac{u^2}{-t}$ and $\frac{(Lu)^2}{-t}$ yields
$$\int_{(B_1\setminus B_r)\times (-r^2,0]}\left\{|\nabla u|^2+|u|+\frac{u^2+(Lu)^2}{-t}\right\}\d\gamma\leq Cr^4.$$
\noindent {\bf Step 4: Conclusion}\\
Combining the estimates of the two contributions (from Steps 2 and 3), the result follows.
\qeda

\section{Monotonicity formulae and classification of free boundary points}\label{sec:mon}

\subsection{A Weiss type monotonicity formula}

Let $u$ be a solution of \eqref{eq:obstacle}. Define the Weiss energy as
\begin{align}\label{eq:Wdef}
W(r,u)=E(r,u)-\int_0^r\frac{2}{s^5}\int_{S_s}Lu(1-Hu)\d\gamma \d s,
\end{align}
with 
$$E(r,u)=\frac{1}{r^4}\int_{S_r}\left(|\nabla u|^2+2u+\frac{u^2}{t}\right)\d\gamma.$$
Then we have the following result:

\begin{prop}{\bf (Weiss type Monotonicity formula)}\label{prop:weiss}\\
Let $u$ be a solution of \eqref{eq:obstacle} for $p>(n+2)/2$ with $p\ge 2$ such that $\tilde{\sigma}_p(1)$ is finite.
Then $E(r,u)$ is well defined for any $r\in (0,1]$.
Furthermore, if $\sigma_p$ is Dini, then $W(r,u)$ is also well defined for any $r\in (0,1]$.
Moreover $W(r,u)$ is non-decreasing in $r$ and we have for any $0<s<\tau\le 1$:
\begin{equation}\label{eq::el36}
W(\tau,u)-W(s,u)=\int_{s}^\tau \frac{dr}{r^5}\int_{S_r}(Lu)^2\d\gammab.
\end{equation}
\end{prop}

\begin{rem} 
In the case where $f$ is constant equal to $1$, and $u$ is a solution of \eqref{eq:obstacle} 
with polynomial growth at infinity then $W=E$ is non-decreasing and constant if and only if 
$u$ is parabolically homogeneous of degree 2. This is proved in Theorem 3.1 in \cite{Wei99}. 
Moreover, a similar formula has also been obtained in \cite{Bla06} and \cite{CPS04}.
\end{rem}

\noindent {\bf Proof of Proposition \ref{prop:weiss}}\\
We split the proof into several steps.\\
\noindent {\bf Step 1: Definition of $E$ and $W$}\\ 
Notice that since $\tilde{\sigma}_p(1)$ is finite, we can apply Corollary \ref{cor:lpprim}.
Using the fact that the support of $u$ is contained in $Q_1^-$, we see that the integral over $S_r$ reduces to an integral over $B_1\times (-r^2,0]$ and we can thus 
deduce that $E$ is well defined and bounded.

For things to make sense we also have to prove that $W$ is well defined,
i.e., that all the involved terms are finite. We start by proving that
the integrand in the second term of \eqref{eq:Wdef} is finite.
Since $u$ has support contained in $B_1\times (-1,0]$, $Hu=f\chi_{\left\{u>0\right\}}$ and $f(0)=1$ it reads
\begin{align*}
 \frac{1}{s^5} \int_{B_1\times (-s^2,0]}Lu(f(0)-f)\d\gamma=A(s).
\end{align*}
We fisrt notice that $A(r)$ is well defined, because $f\in L^p$ and $Lu\in L^{p'}$ by Corollary \ref{cor:lpprim}.
More quantitatively, we get from Corollary \ref{cor:lpprim} and H\"older's inequality 
\begin{align*}\nonumber 
 &\int_0^r |A(s)|\d s \nonumber \\
 &\leq\int_0^r s^{-3}\left(\frac{1}{s^2}\int_{B_1\times (-s^2,0]}
   |f-f(0)|^p\d\gamma \right)^\frac{1}{p}\left(\frac{1}{s^2}\int_{B_1\times (-s^2,0]}
     |Lu|^{p'}\d\gamma\right)^\frac{1}{p'}\d s\\
&\leq C \int_0^r \frac{\sigma_p(s)}{s} \d s<\infty,\nonumber
\end{align*}
if $\sigma_p$ is Dini.\\
\noindent {\bf Step 2: Formal computations}\\ 
We observe that
$$E(rs,u)=E(s,u_r),$$
with 
$$u_r(x,t)=\frac{u(rx,r^2t)}{r^2}.$$
Moreover, we introduce the notation
$$Lu:=\frac{\d}{\d r}u_r\Big|_{r=1}=x\cdot \nabla u+2tu_t-2u.$$
Using that $\nabla G(x,-t)=\frac{x}{2t}G(x,-t)$, integration by parts gives formally
\begin{align}
rE'(r,u)&=\frac{\d}{\d
  s}E(sr,u)\Big|_{s=1}\nonumber\\
&=\frac{\d}{\d
  s}E(r,u_s)\Big|_{s=1}\nonumber\\
&\nonumber=\frac{1}{r^4}\int_{S_r}\left(2\nabla
u\cdot \nabla Lu+2Lu+\frac{2uLu}{t}\right)G(x,-t)\d x \d t \\
&=\frac{1}{r^4}\int_{S_r}Lu(-2(Hu-1)-\frac{Lu}{t})G(x,-t)\d
x\d t\label{eq:Eprim}
\\
&=\nonumber\frac{2}{r^4}\int_{S_r}Lu(1-Hu)G(x,-t)\d
x \d t+\frac{1}{r^4}\int_{S_r}\frac{(Lu)^2}{-t}G(x,-t)\d
x\d t.
\end{align}
\noindent {\bf Step 3: Justifying the formal computations}\\
Our goal is now to justify the formal computations done in Step 2.
Notice the presence of the term $\nabla Lu$ in the integral of the third line of (\ref{eq:Eprim}),
where $Lu$ contains terms like $u_t$ whose best regularity is  $L^p$. This shows that the third line of (\ref{eq:Eprim}) has no meaning. Nevertheless after a formal 
integration by parts, the fourth line of (\ref{eq:Eprim}) makes perfect sense.

We now propose a way to justify the formal computations above, which consists in replacing $S_r$ everywhere by
$$S^\delta_r=\R^n\times (-r^2,-r^2\delta)$$
for some fixed $\delta\in (0,1)$ (that will tend to zero at the end of the reasoning). We then denote by $E^\delta$ and $W^\delta$, the associated quantities
where in $E$ and $W$, we replace $S_r$ by $S_r^\delta$. Notice that our choice $S^\delta_r$ still respects the desired scale invariance.
Now consider a mollifier $\rho_\varepsilon(x,t)$ with support in $Q_\varepsilon$ and define
$$u_\varepsilon= \rho_\varepsilon\star u \quad \mbox{on}\quad \R^n\times (-\infty,-\varepsilon^2)$$
Then the computation for $E^\delta$ similar to (\ref{eq:Eprim}) gives without any difficulties for $r^2 \delta >\varepsilon^2$
$$\left(E^\delta\right) '(r,u_\varepsilon) = \frac{2}{r^5}\int_{S_r^\delta}Lu_\varepsilon(1-Hu_\varepsilon)G(x,-t)\d
x \d t+\frac{1}{r^5}\int_{S_r^\delta}\frac{(Lu_\varepsilon)^2}{-t}G(x,-t)\d x\d t.$$
and then for all $0<\varepsilon \delta^{-\frac12}<s<\tau\le 1$
$$E^\delta(\tau,u_\varepsilon)-E^\delta(s,u_\varepsilon) = \int_s^\tau  \left(\frac{2}{r^5}\int_{S_r^\delta}Lu_\varepsilon(1-Hu_\varepsilon)G(x,-t)\d
x \d t+\frac{1}{r^5}\int_{S_r^\delta}\frac{(Lu_\varepsilon)^2}{-t}G(x,-t)\d x\d t\right)\d r.$$
Using the convergence $u_\varepsilon\to u$ in $W^{2,1}_p(\R^n\times (-\infty,-\kappa))$ for any $\kappa>0$,
we get for all $0<s<\tau\le 1$:
\begin{equation}\label{eq::el35}
E^\delta(\tau,u)-E^\delta(s,u) = \int_s^\tau  \left( \frac{2}{r^5}\int_{S_r^\delta}Lu(1-Hu)G(x,-t)\d
x \d t+\frac{1}{r^5}\int_{S_r^\delta}\frac{(Lu)^2}{-t}G(x,-t)\d x\d t\right) \d r.
\end{equation}
We can then simply let $\delta$ tend to zero and pass to the limit in (\ref{eq::el35}), using the integrability of the integrand defining $E$,
the integrability of the integrand in the first integral of the right hand side of (\ref{eq::el35}),
and the monotone convergence theorem for passing to the limit in the last integral of (\ref{eq::el35}).
This shows (\ref{eq::el36}).
\qeda\\

\subsection{Classification of blow-up limits and free boundary points}

In this section we give the proof of Theorem \ref{thm:class}. 
The proof is standard once we have the monotonicity formula at our hands, 
but for the matter of clarity we present a brief version.\\

\noindent {\bf Proof of Theorem \ref{thm:class}}\\
Let $u$ be a solution of \eqref{eq:obstacle} and put
$$u_r(x,t)=\frac{u(rx,r^2t)}{r^2}.$$
By Proposition \ref{prop:quadgrowth}, Theorem \ref{thm:globalpar} and
Theorem \ref{thm:sobolev}, the functions
 $u_r$ are locally uniformly bounded in
$W^{2,1}_p(\R^n\times \R^-)\cap C^\alpha(\R^n\times \R^-)$ for some $\alpha>0$ and for all
$r\in (0,1]$. Thus, we can extract a subsequence $u_{r_j}$ with $r_j\to 0$, converging locally uniformly to $u_0$ where the limit $u_0$ satisfies
(using $W^{2,1}_p$ estimates to check that $Hu_0=0$ a.e. on $\left\{u_0=0\right\}$)
$$Hu_0=f(0)\chi_{\{u_0>0\}} \text{ in }\R^n\times \R^-$$
and $u_0$ has at most quadratic growth at infinity (as a consequence of Proposition \ref{prop:quadgrowth}). 
From Proposition \ref{prop:weiss} it follows that the limit
$$W(0^+,u):=\lim_{r\searrow 0}W(r,u)=\lim_{r\searrow 0}E(r,u)$$ 
exists. In addition, for all $S>0$
\begin{align*}
W(0^+,u)&=\lim_{j\to\infty}W(r_j,u)=\lim_{j\to\infty} W(Sr_j,u).
\end{align*}
Therefore, by Proposition \ref{prop:weiss}
(and the lower semi-continuity property of the norm  for the weak convergence in the Hilbert space of functions $h$ 
such that the square of the norm of $h$ is given by  $\int_R^S \frac{1}{r^5}\int_{S_r}h^2d\gammab dr < +\infty$)
\begin{align*}
0&=\lim_{j\to\infty} \left(W(S,u_{r_j})-W(R,u_{r_j})\right)\\&=\lim_{j\to\infty}\int_R^S
\frac{1}{r^5}\int_{S_r}(Lu_{r_j})^2\d\gammab \d r\geq \int_R^S
\frac{1}{r^5}\int_{S_r}(Lu_{0})^2\d\gammab \d r,
\end{align*}
whenever $S>R>0$. Hence, $Lu_0=0$ and it follows that $u_0$ must be parabolically homogeneous of degree 2. 
By Lemma 6.2 and 6.3 in \cite{CPS04}, $E(u_0)$ can only take the values $0,15/2$ and $15$ and the result follows.
\qeda

\subsection{A monotonicity formula for the singular free boundary points}

Define $R(r,u)$ so that $W(r,u)=E(r,u)-2R(r,u)$, i.e.
\begin{equation}\label{eq::el9}
R(r,u)=\int_0^r\frac{1}{s^5}\int_{S_s}Lu(1-Hu)\d\gamma\d s,
\end{equation}
and let for $r\in (0,1]$
\begin{equation}\label{eq::el7}
\n (r,u)=\frac{1}{r^4}\int_{S_r}w^2\d\gammab,
\end{equation}
with $w=u-v_0$ and $v_0\in \mathcal{P}_\textup{sing}$ (see (\ref{eq::el1})).
Then we have the following result:

\begin{prop}\label{prop:sing}{\bf (Monotonicity formula for singular points)}\\ 
Let $u$ be a solution of \eqref{eq:obstacle} for $p> (n+2)/2$ with $p\ge 2$ such that $\tilde{\sigma}_p(1)$ is finite.
Let $v_0\in \mathcal{P}_\textup{sing}$. Then $\n(r,u)$ is well defined for all $r\in (0,1]$. Assume moreover that $\sigma_p$ is Dini and that the origin is a singular point as defined in Theorem \ref{thm:class}.
Then there exists a function $F\in L^1((0,1))$ such that for any $0<s<\tau\le 1$:
\begin{equation}\label{eq::el43}
\n(\tau,u)-\n(s,u)\geq \int_s^\tau \left(F(r)+\frac2r\int_0^r\frac{da}{a^5}\int_{S_a}(Lu)^2\d\gammab\right) \d r,
\end{equation}  
where in particular for all $\tau\in (0,1]$ we have the estimate
\begin{align*}
\int_0^\tau |F(r)|\ dr\leq 4\int_0^\tau \left( \frac{\sigma_2(r)}{r}\left(\frac{1}{r^6}\int_{B_1\times
  (-r^2,0]}|u-v_0|^{2}\d \gamma\right)^\frac{1}{2}\right)\d r +32\left(\int_0^\tau \frac{\sigma_2(r)}{r}\ dr\right)^2.
\end{align*}
Moreover we have for all $r\in (0,1]$
\begin{equation}\label{eq::el11}
\left(\frac{1}{r^6}\int_{B_1\times
  (-r^2,0]}|u-v_0|^{2}\d \gamma\right)^\frac{1}{2} \le C_1
\end{equation}
for some constant $C_1=C_1(||f||_{L^p(Q_1^-)},||u||_{L^p(Q_1^-)},\tilde{\sigma}_p(1),n,p)$.
\end{prop}

\begin{rem} 
This sort of monotonicity formula specially designed for singular points 
was introduced in \cite{Mon03} for the elliptic obstacle problem. 
For the parabolic case, a similar monotonicity formula was proved in \cite{Bla06}. 
The major difference is that in the formula in \cite{Bla06}
there is no integration in time, contrarily to our formula (\ref{eq::el7}).
\end{rem}

\noindent {\bf Proof of Proposition \ref{prop:sing}}\\
We first notice that $\n$ is well defined because $w$ has a quadratic growth, i.e. satisfies estimate (\ref{eq::el6}) 
and then we can apply Corollary \ref{cor:lpprim}.
We also observe for future use that
\begin{equation}\label{eq::el8}
\n(rs,u)=\n(s,u_r),
\end{equation}
where as before
$$u_r(x,t)=\frac{u(rx,r^2t)}{r^2}.$$
Since the rest of the proof is quite long, we split into several steps.\\
\noindent {\bf Step 1: Formal computations}\\ 
Using the scaling property (\ref{eq::el8}) of $\n$, we compute {\it formally}
\begin{equation}\label{eq::el9bis}
r \n'(r,u)=\frac{\d}{\d s}\n(rs,u)\Big|_{s=1}=\frac{\d}{\d s}\n(r,u_s)\Big|_{s=1}=\frac{2}{r^4}\int_{S_r} \frac{w(Lw)}{-t}G(x,-t)\d x\d t.
\end{equation}
With $E$ as before we have {\it formally} by integration by parts
\begin{align}\nonumber 
E(r,u)-E(r,v_0)& \displaystyle =\frac{1}{r^4}\int_{S_r}\left( \nabla(u+v_0)\cdot \nabla w+2w+\frac{(u+v_0)}{t}w\right)G\ \d x\d t
\\\label{eq::el38}
&\displaystyle =\frac{1}{r^4}\int_{S_r}\left (-\lap
  (u+v_0)-\frac{x}{2t}\cdot \nabla (u+v_0)+2+\frac{(u+v_0)}{t}\right)wG \ \d x\d t
\\\nonumber
&\displaystyle =\frac{1}{r^4}\int_{S_r}(-Hu+1-\frac{Lu}{2t})w G \ \d x\d t,
\end{align}
where we used that $Lv_0=0$ since $v_0$ is parabolically homogeneous of degree 2 and that $Hv_0=1$. 
Together with (\ref{eq::el9}), this implies
\begin{equation}\label{eq::el40}
\begin{array}{ll}
\n'(r,u) &\displaystyle =4\left(\frac{E(r,u)-E(r,v_0)}{r}\right)+\frac{4}{r^5}\int_{S_r}(Hu-1)w\d\gamma\\
\\
&\displaystyle =4\left(\frac{E(r,u)-E(r,v_0)}{r}\right)+4Q(r,u)+ 4Q_0(r,u),
\end{array}\end{equation}
with
$$Q(r,u)=\frac{1}{s^5}\int_{S_s\cap \left\{u>0\right\}}(Hu-1)w  \d\gamma \quad \mbox{and}\quad 
Q_0(r,u)=\frac{1}{s^5}\int_{S_s\cap \left\{u=0\right\}}v_0\d\gamma\ge 0.$$
Notice that $Q(r,u)$ and $Q_0(r,u)$ are well defined, as a consequence of Corollary \ref{cor:lpprim} 
and H\"{o}lder's inequality as in Step 1 of the proof of Proposition 
\ref{prop:weiss}.
Using the definition (\ref{eq::el9}) of $R(r,u)$, we get
\begin{align*}
\frac{E(r,u)-E(0,u)}{r}&=\frac{W(r,u)-W(0,u)}{r}+\frac{2R(r,u)}{r}\\
&=\frac{2R(r,u)}{r}+\frac{I(r,u)}{r}.
\end{align*}
where we have used Proposition \ref{prop:weiss} for the last line
with the finite quantity
\begin{equation}\label{eq::el14}
I(r,u)=\int_0^r\frac{1}{s^5}\int_{S_s}(Lu)^2\d\gammab \d s.
\end{equation}
Since $E(0^+,u)= E(r,v_0)$ for all $r>0$ 
(by the homogeneity of $v_0$ and using the fact that the origin is a singular point as in Theorem \ref{thm:class}), we deduce
\begin{equation}\label{eq::el41}
\frac{E(r,u)-E(r,v_0)}{r} = \frac{2R(r,u)}{r}+\frac{I(r,u)}{r}.
\end{equation}
Therefore
\begin{equation}\label{eq::el10}
\n'(r,u)=4Q(r,u)+4Q_0(r,u)+\frac{8R(r,u)}{r}+\frac{4I(r,u)}{r}.
\end{equation}

\noindent {\bf Step 2: Estimates}\\
We first estimate the contribution from $Q(r,u)$. 
By Corollary \ref{cor:lpprim} and H\"older's inequality
\begin{align*}
\int_0^r|Q(s,u)|\d s&= \int_0^r \frac{1}{s^5}\int_{\R^n\times
  (-s^2,0]}\chi_{\left\{u>0\right\}}|f-f(0)||u-v_0|\d \gamma\d s\\
&\le  \int_0^r \frac{1}{s^5}\int_{B_1\times
  (-s^2,0]}|f-f(0)||u-v_0|\d \gamma\d s\\
&\leq \int_0^r  \frac{1}{s}\left(\frac{1}{s^2}\int_{B_1\times
  (-s^2,0]}|f-f(0)|^2\d \gamma\right)^\frac{1}{2}\left(\frac{1}{s^6}\int_{B_1\times
  (-s^2,0]}|u-v_0|^{2}\d \gamma\right)^\frac{1}{2}\d s
\end{align*}
we obtain
\begin{equation*}
\int_0^r|Q(s,u)|\d s\le  \int_0^r \frac{\sigma_2(s)}{s} \left(\frac{1}{s^6}\int_{B_1\times
  (-s^2,0]}|u-v_0|^{2}\d \gamma\right)^\frac{1}{2}\d s.
\end{equation*}
Using the uniform quadratic growth of  all elements $v_0\in \ps$ together with the first estimate of Corollary \ref{cor:lpprim} for $p'=2$, we obtain (\ref{eq::el11}).

It remains to prove that the sum of the two last terms of (\ref{eq::el10}), i.e.,
$$\frac{8R(r,u) + 4I(r,u)}{r}$$
can be bounded from below by something integrable. As in the proof of Proposition \ref{prop:weiss}, we use the notation $A(s)$:
\begin{align*}
R(r,u)= \int_0^rA(s) \d s \quad \mbox{with}\quad A(s)=\frac{1}{s^5}\int_{S_s}Lu(f(0)-f))\d\gamma.
\end{align*}
The idea is to absorb some bad parts of this term into a part of the other
term $I(r,u)$. We have with $g=f-f(0)$, by H\"older's and Young's inequalities
\begin{align*}
|A(s)| & \leq \frac{1}{s^5}\int_{S_s}
\sqrt{-t}|g|\frac{|Lu|}{\sqrt{-t}}\d\gamma \\
&\leq \left(\frac{1}{s^5}\int_{S_s}|g|^2(-t) \d\gamma\right)^\frac{1}{2}\left(\frac{1}{s^5}\int_{S_s}
  \frac{|Lu|^{2}}{-t} \d\gamma\right)^\frac{1}{2}\\
&\leq \e\left(\frac{1}{s^5}\int_{S_s}|g|^2(-t) \d\gamma
\right)+\frac{1}{\e}\left(\frac{1}{s^5}\int_{S_s}
  |Lu|^{2}\d\gammab\right).
\end{align*}
Therefore, with $\e=4$, we obtain
$$|A(s)|-\frac14\frac{1}{s^5}\int_{S_s}(Lu)^2\d\gammab  \leq 4\left(\frac{1}{s^5}\int_{S_s}|g|^2(-t)\d\gamma\right)\leq 4\left(\frac{1}{s^3}\int_{S_s}|g|^2\d\gamma\right).$$
Thus
\begin{equation}\label{eq::el13}
R(r,u)+\frac{1}{4}I(r,u) \ge -\int_0^r  \left(|A(s)| -\frac14\frac{1}{s^5}\int_{S_s}(Lu)^2\d\gammab\right)\d s \ge -4 \int_0^r  \left(\frac{1}{s^3}\int_{S_s}|g|^2\d\gamma\right)\d s
\end{equation}
and
$$\int_0^\tau\frac{1}{r}\int_0^r \frac{1}{s^3}\int_{S_s}|g|^2\d\gamma \d s\d r\le \int_0^\tau\frac{1}{r}\int_0^r\frac{(\sigma_2(s))^2}{s}\d s\d r\le \left(\int_0^\tau\frac{\sigma_2(s)}{s}\d s\right)^2$$
which is finite if $\sigma_2$ is Dini.

\noindent {\bf Step 3: Formal conclusion}\\ 
Using (\ref{eq::el10}) and combining (\ref{eq::el13}) with the fact that $Q_0(r,u)\ge 0$, we get that
$$\n'(r,u)\geq F(r)+\frac{2 I(r,u)}{r},$$
with $I(r,u)$ given in (\ref{eq::el14}) and  where we can take
$$F(r)=4Q(r,u)-\frac{32}{r}\int_0^r\frac{1}{s^3}\int_{S_s}|g|^2\d\gamma \d s$$
so that
\begin{align*}
\int_0^\tau |F(r)|\ dr
\leq 
4 \int_0^\tau \left( \frac{\sigma_2(s)}{s} \left(\frac{1}{s^6}\int_{B_1\times
  (-s^2,0]}|u-v_0|^{2}\d \gamma\right)^\frac{1}{2}\right) \d s
+32\left(\int_0^\tau \frac{\sigma_2(r)}{r}\d r\right)^2.
\end{align*}

\noindent {\bf Step 4: Making everything rigorous}\\ 
Our goal is to justify the formal computations of Step 1 (especially (\ref{eq::el9bis}) and the integration by parts in (\ref{eq::el38})),
and also the formal reasoning in Step 3.
To this end, we proceed as in Step 3 of the proof of Proposition \ref{prop:weiss}.
In every quantity $\n,E,W,R,I,Q,Q_0,F$, we replace $S_r$ by $S_r^\delta$ for some $\delta>0$, 
and define $\n^\delta,E^\delta,W^\delta,R^\delta,I^\delta,Q^\delta,Q_0^\delta,F^\delta$.

Considering a mollification $u_\varepsilon$ of $u$, we first get (\ref{eq::el38}) at the level $\delta$ and for $u_\varepsilon$.
Indeed we can perform the integration by parts in that case, because we can show that the boundary term does not contribute.
To this end, we have to  use the $L^\infty$ bound on $u_\varepsilon$ and $\nabla u_\varepsilon$ 
as a consequence of the $L^\infty$ bound on $u$ given in (\ref{eq::el6}). The other terms are coming from the polynomial $v_0$.
Considering the integration by parts on a bounded domain, we can easily check that the boundary term (which behaves like a polynomial times a Gaussian) goes to zero, see Lemma \ref{lem:ipp}.

Using (\ref{eq::el38}) at the level $\delta$ for $u_\varepsilon$, we deduce the analogue of (\ref{eq::el40}), namely for all $0< s<\tau\le 1$:
$$\n^\delta(\tau,u)-\n^\delta(s,u) = \int_s^\tau \left(4\left(\frac{E^\delta(r,u)-E^\delta(r,v_0)}{r}\right)+4Q^\delta(r,u)+ 4Q_0^\delta(r,u)\right)\d r$$
that we get first for $u_\varepsilon$ for all  $0<\varepsilon \delta^{-\frac12} < s<\tau\le 1$, 
and then for $u$ (after passing to the limit $u_\varepsilon\to u$).
Using the analogue of (\ref{eq::el41}) at the level $\delta$, we deduce the integral form of  (\ref{eq::el10}), i.e.
$$\n^\delta(\tau,u)-\n^\delta(s,u)=\int_s^\tau  \left(4Q^\delta(r,u)+4Q_0^\delta(r,u)+\frac{8R^\delta(r,u)}{r}+\frac{4I^\delta(r,u)}{r}\right)\d r.$$
Using the analogue of (\ref{eq::el13}) at the level $\delta$, we get
$$4Q^\delta(r,u)+ \frac{8R^\delta(r,u)}{r}+\frac{4I^\delta(r,u)}{r} \ge F^\delta(r) + \frac{2I^\delta(r,u)}{r}.$$
Using moreover the fact that $Q_0^\delta(r,u)\ge 0$, we deduce
\begin{equation}\label{eq::el42}
\n^\delta(\tau,u)-\n^\delta(s,u)\ge \int_s^\tau  \left(F^\delta(r) + \frac{2I^\delta(r,u)}{r}\right)\d r.
\end{equation}
We then let $\delta$ tend to zero and pass to the limit in (\ref{eq::el42}), using the integrability of the integrand defining $\n$ and $F$,
and using the monotone convergence theorem for the integral of the last term in (\ref{eq::el42}).
This gives (\ref{eq::el43}).
\qeda\\

As a corollary, there is a unique blow-up at singular free boundary points.
\begin{cor}\label{cor:unique} {\bf (Uniqueness of blow-ups at singular free boundary points)}\\ 
Under the same assumptions as in Proposition \ref{prop:sing} there is a unique blow-up limit at the origin,
i.e., $$u_r(x,t)=\displaystyle \frac{u(rx, r^2t)}{r^2} $$ has a unique limit as $r\to 0$.
\end{cor}

\noindent {\bf Proof of Corollary \ref{cor:unique}}\\
Recall that the set of all possible blow-ups at singular points is the set $\P_{\textup{sing}}$ (see Lemma 6.3 in \cite{CPS04}).
We now fix $v_0\in \ps$ such that there exists a subsequence of rescaling $u_{r_j}$  converging to $v_0$,
and use this function $v_0$ for the definition of $\n(r,u)$.
We define 
$$M(r,u)=\n(r,u)-\int_0^rF(s),$$
where $F$ is as in Proposition \ref{prop:sing}. 
By the latter, $M$ is monotone in $r$ so that the limit
$$M(0^+,u)=\lim_{r\searrow 0}M(r,u)$$
exists. Then, using the scaling invariance of $\n$ (similar to the one used in the proof of Theorem \ref{thm:class}), we get
$$M(r_j,u)=\n(r_j,u)-\int_0^{r_j}F=\n(1,u_{r_j})-\int_0^{r_j}F\to \n(1,v_0) = 0.$$
The convergence of $\n(1,u_{r_j})$ to $\n(1,v_0)$ follows from Lebesgue  dominated convergence theorem, because we have the domination
$$|u_{r_j}(x,t)|\le C (|x|^2 + |t|) = C d^2(x,t) \quad  \mbox{with}\quad d\in L^2_{\gammab}(S_1)$$
which in turn follows from the quadratic bound (\ref{eq::el6}).
Therefore, $M(0^+,u)=0$. Suppose now that for another subsequence $s_j\to 0$, $u_{s_j}\to w_0\neq v_0$. 
Then by Fatou's lemma and Proposition \ref{prop:sing}
$$0=M(0^+,u)=\lim_{j\to \infty}M(s_j,u)=\lim_{j\to\infty}\n(s_j,u)=\lim_{j\to\infty}\n(1,u_{s_j})\geq \n(1,w_0),$$
which implies $w_0=v_0$. 
\qeda

\section{Regularity results when $f$ is Dini}\label{sec:dini}
In this section we will make the stronger assumption that $f$ is Dini
continuous, i.e., that there exists a Dini modulus of continuity $\sigma$ such that
$$|f(Y)-f(X)|\le \sigma(d(Y-X)).$$
In particular $f$ is continuous. 
One important
result that holds in this case (when $f>0$) is the fact that the solution does  not
decay too fast around free boundary points, the so-called non-degeneracy
property. 

\begin{prop} \label{prop:nondeg2}{\bf (Non-degeneracy)}\\ 
Let $u$ be a solution of \eqref{eq:obstacle} and assume further that $f$ is continuous with $f(0)=1$. 
Then there exists $C=C(n)>0$ and $r_0\in (0,1/2]$ (with $r_0$ depending on $f$) 
such that for all $(x_0,t_0)\in \left\{u=0\right\}\cap \dd\{u>0\}\cap Q_{r_0}^-$, there holds
$$\sup_{Q_r^-(x_0,t_0)}u\geq Cr^2,$$
whenever $r\in (0,r_0]$.
\end{prop}
The proof of this statement is quite standard. See for instance Lemma 5.1 in \cite{CPS04}. Proposition \ref{prop:nondeg2} implies that there are no degenerate points in the free boundary, 
i.e., that the first case in Theorem \ref{thm:class} does not occur. 
Moreover, it implies that free boundary points are stable in the sense 
that a limit of a sequence of free boundary points of some problems is again a free boundary point of the limit problem.

\subsection{Regular set of the free boundary}

Now we are ready to give the proof of Theorem \ref{thm:mainreg} for regular points, which is standard, 
knowing that the free boundary has no degenerate points, and having a Weiss type monotonicity formula and 
a control of the Taylor expansion at regular points (see Theorem 1.7 in \cite{LM11}).\\

\noindent {\bf Proof of Theorem \ref{thm:mainreg}}\\
\noindent {\bf Step 1: The set $\Gamma_r$ of regular points is open}\\
Since there are no degenerate points, the set of regular points, $\Gamma_r$, 
is the set of points $X_0$ where $W(0^+,\frac{u(X_0+\cdot)}{f(X_0+\cdot)})$ attains its minimum allowed value $15/2$. 
Since $W$ is increasing, the set of regular points is open.\\
\noindent {\bf Step 2: $\Gamma_r$ is locally a graph}\\
Theorem 1.7 in \cite{LM11} implies that there is a unique blow-up limit at each regular free boundary point $X_0=(x_0,t_0)$, i.e.,
\begin{equation}\label{eq:blowupconv}
\frac{1}{f(X_0)}\frac{u(x_0+rx, t_0+r^2t)}{r^2}\longrightarrow  P_{X_0}(x,t) = \frac 12 \left(\max(0,\nu_{X_0}\cdot x)\right)^2,
\end{equation}
locally in $\R^n\times \R^-$ as $r\to 0$, where $\nu_{X_0}$ is a unit vector that can be interpreted as the normal to the free boundary pointing in the direction of $\left\{u>0\right\}$. We remark that the convergence \eqref{eq:blowupconv} remains true also locally in $\R^n\times \R^+$ as long as $t_0<0$, due to the unique continuation property and the fact that the quadratic growth (cf. Proposition \ref{prop:quadgrowth}) can be proved in the whole $Q_r$. 
Moreover from the estimate given in  Theorem 1.7 in \cite{LM11}, we can easily deduce that on the set of regular points, 
the map $X_0\mapsto P_{X_0}$ is continuous
and then the map $X_0\to \nu_{X_0}$ is also continuous.

Upon rotating the coordinates, we can assume that $\nu_{X_0}=e_n$. Moreover, since any blow-up of $u$ at $X_0$ is of the form $P_{X_0}(x,t)=\frac 12 \left(\max(0,x_n)\right)^2$,
we can locally define for $x'=(x_1,...,x_{n-1})$
$$\tilde{f}(x',t)=\inf\left\{x_n,\quad u(x',x_n,t)>0\right\}.$$
Then near $X_0$ 
$$u>0 \quad \mbox{in}\quad \left\{x_n> \tilde{f}(x',t)\right\},$$
and we also claim that near $X_0$
\begin{equation}\label{eq::el16}
u=0 \quad \mbox{in}\quad \left\{x_n\le  \tilde{f}(x',t)\right\}.
\end{equation}
Indeed, if (\ref{eq::el16}) is not true locally, then we can find a sequence of points $ X^k = (x^k,t^k)\in \Gamma_r$
such that $(x^k)_n< \tilde{f}((x^k)',t^k)=:\tilde{x}^k_n$, i.e. the point $X^k$ is below the point $\tilde{X}^k=((x^k)',\tilde{x}^k_n,t^k)$ 
in the direction of $e_n$. We can also assume that $X^k,\tilde{X}^k\to X_0$. 
Since $u$ is uniformly close to its blow-up limit both at $X^k$ and at $\tilde{X}^k$ 
with $\nu_{X^k}, \nu_{\tilde{X}^k}\to \nu_{X_0} = e_n$, we easily get a contradiction for a rescaling of $u$ at the scale $d(X^k-\tilde{X}^k)$.
Therefore the free boundary $\Gamma_r$ is locally the graph of $\tilde{f}$.\\
\noindent {\bf Step 3: $\Gamma_r$ is locally a $C^1$-graph with respect to $d$}\\
We now want to show that the spatial normal to the graph of $\tilde{f}$ at $X=(x',\tilde{f}(x',t),t)$ is $\nu_X$. In other words, if $g(x',t)$ is defined by \begin{equation}\label{eq::el45}
\nu_{X}=\frac{1}{\sqrt{1+ |g(X')|^2}}(g(X'),1),
\end{equation}
then for $X'=(x',t)$ we wish to show that $g(X')$ is the spatial gradient of $\tilde{f}$. Since $\nu_{X}$ is continuous, this would imply that $\tilde f$ is $C^1$ with respect to the parabolic distance on $\R^{n-1}\times\R$ given by $\sqrt{|x'|^2 + |t|}$. It is sufficient to prove that with $H'=(h',h^t)$, $d'(H')=\sqrt{|h'|^2 + |h^t|}$ and $g(X')$ defined by \eqref{eq::el45}
\begin{equation}\label{eq::el17}
|\tilde{f}(X'+H')-\tilde{f}(X')-g(X')\cdot h'| \le d'(H') \ \omega(d'(H'))
\end{equation}
for some modulus of continuity $\omega$.
Assume by contradiction that (\ref{eq::el17}) is false. Then we can find sequences $X^k, Y^k \to X_0$ such that
$$\left\{\begin{array}{l}
X'^k = ((x^k)',t^k),\\
Y'^k=((y^k)',s^k),\\
H'^k=Y'^k-X'^k=(h'^k,(h^t)^k) \to 0,\\
X^k=((x^k)',\tilde{f}((x^k)',t^k), t^k)=(x^k,t^t),\\
Y^k= ((y^k)',\tilde{f}((y^k)',s^k), s^k)=(y^k,s^k),\\
\end{array}\right.$$
and such that
\begin{equation}\label{eq::el18}
\frac{|\tilde{f}(X'^k+H'^k)-\tilde{f}(X'^k)-g(X'^k)\cdot h'^k|}{d'(H'^k)} \ge \delta >0 \quad \mbox{and}\quad H'^k\to 0.
\end{equation}
Since the distance to the blow-up limit at the point $X_0$ decays uniformly (cf. Theorem 1.7 in \cite{LM11}), we have with $r_k = d(Y^k-X^k)$
$$u^k(x,t)=\frac{u(x^k + r_k x,t^k + r_k^2 t)}{r_k^2} \to u^\infty(x,t)= P_{X_0}(x,t)=f(X_0) \frac12 (\max(0,x_n))^2.$$
Moreover, using the non-degeneracy (Proposition \ref{prop:nondeg2}) to conclude that limits of free boundary points are again free boundary points, we have upon choosing a subsequence
\begin{equation}\label{eq::el46}
Z^k=\left(\frac{y^k-x^k}{r_k},\frac{s^k-t^k}{r_k^2}\right) \to Z=(z,\tau)\quad 
\mbox{with}\quad Z\in \partial \left\{u^\infty>0\right\} \quad \mbox{and}\quad d(Z)=1.
\end{equation}
This implies
\begin{equation}\label{eq::el19}
z_n=\nu_{X_0}\cdot z=0 \quad \mbox{with}\quad z=(z',z_n).
\end{equation}
Passing to the limit in (\ref{eq::el18}), we also get
$$|z_n - g(X_0)\cdot z'|\ge \delta d'(Z')\quad \mbox{with}\quad Z'=(z',\tau),$$
where $g(X_0)=0$ because of (\ref{eq::el45}) with $\nu_{X_0}=e_n$.
This implies (using (\ref{eq::el19})) that $z_n=0$ and $Z'=0$, which is in contradiction with (\ref{eq::el46}).
As a conclusion, (\ref{eq::el17}) holds in a neighborhood of $X_0$ and this ends the proof.
\qeda\\

\subsection{Singular set of the free boundary}

Now we turn our attention to the singular part of the free boundary. 
As an improvement of the uniqueness of the blow-up limit (cf. Corollary \ref{cor:unique}), we have:

\begin{cor}\label{cor:coeffcont}{\bf (Continuity of the blow-up limit at singular points)}\\
Let $u$ be a solution of \eqref{eq:obstacle} and suppose that $f$ is Dini continuous with $f>0$.
For any singular point $X_0=(x_0,t_0)$ of the free boundary, denote by $P_{X_0}\in \P_{\textup{sing}}$ the unique blow-up at  $(x_0,t_0)\in \R^n\times \R^-$, i.e., 
of
$$\frac{1}{f(X_0)}\frac{u(x_0+rx,t_0+r^2t)}{r^2} \longrightarrow  P_{X_0}(x,t)=\frac12 {}^t x\cdot Q_{X_0} \cdot x + m_{X_0}t,$$
as $r\to 0$. Then the mapping $X_0\mapsto (Q_{X_0},m_{X_0})$ is continuous on $\Gamma_s$, the set of singular points of the free boundary.
\end{cor}

\noindent {\bf Proof of Corollary \ref{cor:coeffcont}}\\
This is just an adaptation of the proof of Corollary \ref{cor:unique}.
Let us consider a sequence of points $X^k\to X_0$ and we want to show that $P_{X^k}\to P_{X_0}$.
To this end, we simply choose $v_0=P_{X_0}$ in the definition (\ref{eq::el7}) of $\n$, and deduce from Proposition \ref{prop:sing} that there exists a function $F_k$ such that for $0\le \rho\le r$
$$\n(\rho,u(X^k +\cdot))  + \int_\rho^r F_k(s)ds \le \n(r, u(X^k +\cdot))$$
with
$$\int_0^r |F_k(s)|\ ds \le \omega(r),$$
where $\omega$ is independent of $k$. For any $\varepsilon>0$, we choose $r>0$ small enough such that 
$$\omega(r)\le \varepsilon \quad \mbox{and}\quad \n(r, u(X_0 + \cdot))\le \varepsilon,$$
where the last inequality follows from our choice of $v_0$ such that $\n(0^+,u(X_0+\cdot))=0$.
Now for $r$ fixed, we choose $k$ large enough such that
$$\n(r, u(X^k +\cdot)) \le \n(r, u(X_0 + \cdot)) + \varepsilon,$$
which follows from the Lebesgue dominated convergence theorem as in the proof of Corollary \ref{cor:unique}, together with the bound (uniform in $k$)
$$|u(X^k +(x,t))|\le C (|x|^2 + |t|) = C d^2(x,t), \quad d\in L^2_{\gammab}(S_1).$$
This implies that
\begin{equation}\label{eq::el20}
\n(\rho,u(X^k +\cdot))\le 3\varepsilon,
\end{equation}
 for all $0\le \rho\le r$. In particular for $\rho=0$, we obtain
$$\int_{S_1}(P_{X^k}-P_{X_0})^2\ d\gammab  = \n(0^+,u(X^k +\cdot))\le  3\varepsilon.$$
Since $\varepsilon$ is arbitrary, we conclude that $P_{X^k}\to P_{X_0}$, and as a consequence $(Q_{X^k},m_{X^k})\to (Q_{X_0},m_{X_0})$.
\qeda\\

Exploiting the properties further we are able to obtain a finer analysis
of the singular set. For this purpose we need the definition below and
the results that follow.\\

For $\delta>0$, let us denote by $\Gamma_\delta(j+1)$ for $0\le j\le n-1$, the set of points $X$ in $\Gamma_s$ for which the
associated matrix $Q_X$ has at least $n-j$ eigenvalues bounded from below by $\delta$, 
i.e., if we by $\lambda_1(Q_X)\le \lambda_2(Q_X)\le ... \le \lambda_n(Q_X)$ denote the eigenvalues of the matrix $Q_X$ in nondecreasing order, then
$$\Gamma_\delta(j+1)=\left\{X\in \Gamma_s: 0<\delta \le \lambda_{j+1}(Q_X)\le ... \le \lambda_n(Q_X)\right\}.$$
Furthermore, for $X_0\in \Gamma_\delta(j+1)$, we introduce the notation
\begin{align}\label{eq::el33}
&K_j(X_0)=\left\{ \text{the vector space associated to the eigenvalues } \lambda_{1}(Q_{X_0}), ...,\lambda_j(Q_{X_0})\right\},\\
&I_j(X_0)=\left\{  \text{the vector space associated to the eigenvalues }  \lambda_{j+1}(Q_{X_0}), ...,\lambda_{n}(Q_{X_0})\right\}.\nonumber 
\end{align}
where  in each case we consider the sum of the corresponding eigenspaces.
In the particular case where $Q_{X_0}$ has a kernel of dimension $j$, we have $K_j(X_0) = \mbox{Ker}\ Q_{X_0}$ and $I_j(X_0)=\mbox{Im}\ Q_{X_0}$.
For a general point $X_0\in \Gamma_\delta(j+1)$,
we define $p^j_{X_0}=p_{X_0}$ as the orthogonal projection on to $K_j(X_0)$, i.e.
$$p_{X_0}(x,t)=(\mbox{Proj}^\perp_{|K_j(X_0)}(x), t)$$
Then the following result holds true.

\begin{lem}\label{lem:projreg}{\bf (Approximation by the tangent space)}\\
Let $u$ be a solution of \eqref{eq:obstacle} and suppose that $f$ is Dini continuous with $f>0$.
Then for fixed $\delta>0$ and $j\in \left\{0,...,n-1\right\}$, we have
$$\lim_{\e \to 0}\sup_{\stackrel{X,Y\in \Gamma_\delta(j+1)}{d(Y-X)\leq \e}}\frac{d\left(Y-X - p^j_X(Y-X)\right)}{d(Y-X)} = 0.$$
\end{lem}

\noindent {\bf Proof of Lemma \ref{lem:projreg}}\\
The proof is basically the same as the one of Lemma 5.2 in \cite{Bla06}, and also similar to Step 7 of the proof of Theorem \ref{thm:mainreg}.
If the result is false for sequences $X^k,Y^k$ converging to some point $X_0\in \Gamma_\delta(j+1)$, 
we simply consider a blow-up along the sequence $X^k$ at the scale $r_k=d(Y^k-X^k)\to 0$.
Notice that 
$$\left|\frac{|Y-X - p^j_X(Y-X)|}{d(Y-X)} - \frac{|Y-X - p^j_Y(Y-X)|}{d(Y-X)}\right|\le |p^j_X-p^j_Y| \le \omega(d(X-Y))$$
for some modulus of continuity $\omega$. Therefore, upon interchanging the roles of $X^k=(x^k,t^k)$ and $Y^k=(y^k,s^k)$, we can assume that $s^k\le t^k$.
As in the proof of Corollary \ref{cor:coeffcont}, we have  (\ref{eq::el20}). Namely, for any $\varepsilon>0$, 
there exists $r_\varepsilon>0$ and an integer $k_\varepsilon$ such that for any $0\le \rho\le r_\varepsilon$ and any integer $k\ge k_\varepsilon$, 
we have (for $\n$ defined with $v_0=P_{X_0}$)
$$\n(\rho,u(X^k +\cdot))\le 3\varepsilon.$$
Now for $k$ large enough, we have $r_k\le r_\varepsilon$
and we can conclude that
$$u^k(x,t)=\frac{1}{f(X^k)}\frac{u(x^k + r_k x, t^k+r_k^2 t)}{r_k^2} \to P_{X_0}(x,t).$$
Our choice $s^k\le t^k$ allows us to arrive at a contradiction as in the proof of Theorem \ref{thm:mainreg} due to the explicit form of $P_{X_0}$ for past times, i.e.,
$$P_{X_0}(x,t)=\frac12 {}^t x\cdot Q_{X_0}\cdot x + m_{X_0} t \quad \mbox{for}\quad (x,t)\in\R^n\times \R^-.$$
\qeda\\

As a corollary we can prove that, given a singular point $X_0\in \Gamma_\delta(j+1)$, 
there exists a neighborhood of $X_0$, such that 
inside this neighborhood, the set $\Gamma_\delta(j+1)$ is contained in a $(j+1)$-dimensional $C^1$-manifold.

\begin{cor}{\bf ($\Gamma_\delta(j+1)$ is locally contained in a $(j+1)$-dimensional $C^1$-manifold)}\label{cor:graph}\\ 
Let $u$ be a solution of \eqref{eq:obstacle} and suppose that $f$ is Dini continuous with $f>0$.
Then for every singular point $X_0\in \Gamma_\delta(j+1)$ for some $0\le j\le n-1$, there exists $r>0$ such that the set
$$\Gamma_\delta(j+1)\cap Q_r(X_0)$$
is contained in a $(j+1)$-dimensional $C^1$-manifold in $\R^{n+1}$ (in the sense of Definition \ref{defi::C1manif}).
\end{cor}

\noindent {\bf Proof of Corollary \ref{cor:graph}}\\ 
\noindent {\bf Step 1: $\Gamma_\delta(j+1)$ is locally a graph over the subset $P_{X_0}(\Gamma_\delta(j+1))$ of $K_j(X_0)\times \R$}\\ 
We first notice that from Lemma \ref{lem:projreg} it follows that there exists a nondecreasing modulus of continuity $\omega$ such that 
with $p_X=p^j_X$ there holds
\begin{equation}\label{eq::el21}
d(Y-X - p_X(Y-X)) \le d(Y-X) \omega(d(Y-X)),
\end{equation}
for all $X,Y\in  \Gamma_\delta(j+1)$. Moreover, from the continuity of the map $X\mapsto Q_X$ on the singular set (cf. Corollary \ref{cor:coeffcont})
\begin{equation}\label{eq::el22}
d((p_Y-p_X)(Z))\le  d(Z)\omega(d(Y-X)) \quad \mbox{for all}\quad Z.
\end{equation}
In particular this implies the continuity of the projection matrices on $\Gamma_\delta(j+1)$, i.e., the continuity of the map $X\mapsto p_X$ on $\Gamma_\delta(j+1)$.
From (\ref{eq::el22}), we deduce that on $\Gamma_\delta(j+1)$
$$d(p_{X}(Y-X)-p_{X_0}(Y-X)) \le d(Y-X) \omega(d(X-X_0)),$$
which implies 
$$d((Y-X)-p_{X_0}(Y-X))\le d(Y-X) \left\{\omega(d(Y-X))+ \omega(d(X-X_0))\right\},$$
where we have used (\ref{eq::el21}). From the triangle inequality
$$d(Y-X)\le d((Y-X)-p_{X_0}(Y-X)) + d(p_{X_0}(Y-X))$$
we deduce
$$ d(Y-X)\left(1- \left[\omega(d(Y-X))+ \omega(d(X-X_0))\right]\right) \le d(p_{X_0}(Y-X)).$$
Now choose $r_0$ such that $\omega(2\sqrt{2}r_0)\le 1/4$. Then if $X,Y\in Q_{r_0}(X_0)$,
we have $d(Y-X_0), d(X-X_0)\le \sqrt{2}r_0$, and thus
\begin{equation}\label{eq::el26}
d(p_{X_0}(Y-X)) \ge \frac12 d(Y-X).
\end{equation}
In particular this implies that $\Gamma_\delta(j+1)\cap Q_{r_0}(X_0)$ is a graph on $p_{X_0}(\Gamma_\delta(j+1))$.\\
\noindent {\bf Step 2: Definition and properties of the map $F$}\\
Considering the decomposition (see (\ref{eq::el33}))
\begin{equation}\label{eq::el25}
\R^{n+1} = (K_j(X_0)\times \R) \oplus I_j(X_0),
\end{equation}
we now set 
$$\left\{\begin{array}{l}
\bar X = p_{X_0}(X),\quad X=(\bar X, F(\bar X)),\\
\bar Y = p_{X_0}(Y),\quad Y=(\bar Y, F(\bar Y))
\end{array}\right.$$
such that $\Gamma_\delta(j+1)\cap Q_{r_0}(X_0)$ is the graph of the function $F:E\to I_j(X_0)$ (from Step 1) defined on 
$$E=\mbox{Proj}^\perp_{|K_j(X_0))\times \R }\left(\Gamma_\delta(j+1)\cap Q_{r_0}(X_0)\right)$$
Using the decomposition (\ref{eq::el25}), we can write the $n\times n$ matrix $p_X$ (associated to the projection also denoted by $p_X$):
$$p_X= \left(\begin{array}{ll}
A_X & B_X\\
{}^tB_X & C_X
\end{array}\right)$$
We deduce from (\ref{eq::el22}) that
\begin{equation}\label{eq::el27}
|A_Y-A_X|, |B_Y-B_X|, |C_Y-C_X|\le   d(Y-X) \omega(d(Y-X))
\end{equation}
In (\ref{eq::el21}), taking the projection on $I_j(X_0)$ of the argument of the term in the left hand side, we deduce that
$$\left|F(\bar Y) - F(\bar X) - ({}^t B_X\ C_X)\cdot \left(\begin{array}{l}
\bar Y - \bar X\\
F(\bar Y)- F(\bar X)
\end{array}\right)\right| \le d(Y-X) \omega(d(Y-X))$$
which implies
$$|F(\bar Y) - F(\bar X) -(I-C_X)^{-1}\ {}^tB_X(\bar Y-\bar X)|\le |(I-C_X)^{-1}|d(Y-X) \omega(d(Y-X))$$
i.e. using (\ref{eq::el26})
\begin{equation}\label{eq::el28}
|F(\bar Y) - F(\bar X) -G_{\bar X} (\bar Y-\bar X)| \le 2 M_0 d(\bar Y-\bar X) \omega(2d(\bar Y-\bar X)),
\end{equation}
where
$$
 G_{\bar X} = (I-C_X)^{-1}\ {}^tB_X
$$
and where we have used the following estimate
$$|(I-C_X)^{-1}| \le M_0 \quad \mbox{for all}\quad X\in Q_{r_0}(X_0),$$
which is valid upon reducing $r_0$ and using (\ref{eq::el27}) with $C_{X_0}=0$.
In particular, the map $\bar X \to G_{\bar X}$ is continuous on $E$ (by composition of applications).\\
\noindent {\bf Step 3: Conclusion}\\
From (\ref{eq::el28}) and the parabolic $C^1$ extension of Whitney type (Proposition \ref{prop:whitney}), 
we conclude that, possibly for yet a smaller $r_0$,  the map $\bar X\mapsto \xi\cdot F(\bar X)$ 
defined on $E$ admits a $C^1$ extension $\tilde{F}_\xi: B_{r_0}(p_{X_0}(X_0))\to \R$, where 
$B_{r_0}(p_{X_0}(X_0)) \subset K_j(X_0)\times \R$.
Let now $\xi_1,...,\xi_{n-j}$ be an orthonormal basis of $I_j(X_0)$ and let
$$\begin{array}{llll}
\tilde{F}: & B_{r_0}(p_{X_0}(X_0)) & \longrightarrow & I_j(X_0)\\
& \bar X & \longmapsto & \displaystyle \tilde{F} = \sum_{i=1,...,n-j} \tilde{F}_{\xi_i}(\bar X)\xi_i
\end{array}$$
We see that the graph of the map $\tilde{F}$ defines a $(j+1)$-dimensional $C^1$-manifold in $\R^{n+1}$ in the sense of Definition \ref{defi::C1manif}.
\qeda\\

For $\delta>0$, let us denote by $\Gamma_\delta(n+0)$, the set of points $X$ in $\Gamma_s$ for which the
associated coefficient $m_X$ (defined in Corollary \ref{cor:coeffcont}) is bounded from above by $-\delta$, i.e.
$$\Gamma_\delta(n+0)=\left\{X\in \Gamma_s: -m_X\ge \delta >0\right\}.$$

\begin{prop}{\bf ($\Gamma_\delta(n+0)$ is locally contained in a $C^2$-graph)}\label{prop:graphn}\\
Let $u$ be a solution of \eqref{eq:obstacle} and suppose that $f$ is Dini continuous with $f>0$.
Then for every singular point $X_1\in \Gamma_\delta(n+0)$, there exists $r>0$ such that the set
$$\Gamma_\delta(n+0)\cap Q_r(X_1)$$
is contained in a (euclidean) $C^2$-graph of the form $t=t(x_1,\ldots,x_n)$.
\end{prop}

\noindent {\bf Proof of Proposition \ref{prop:graphn}}\\
\noindent {\bf Step 1: First estimate}\\
We claim that we have the following estimate
\begin{equation}\label{eq::el30}
\lim_{\e \to 0}\sup_{\stackrel{(x,t),(y,s)\in \Gamma_\delta(n+0)}{d((y,s)-(x,t))\leq \e}}\frac{|s-t|}{|y-x|^2} = 0.
\end{equation}
The proof is similar to the proof of Lemma \ref{lem:projreg} and also similar to Step 7 of the proof of Theorem \ref{thm:mainreg}. We assume that (\ref{eq::el30}) is false for sequences $X^k=(x^k,t^k),Y^k=(y^k,s^k)\to X_0$, where by symmetry, we can always assume that $s^k\le t^k$.
Then we consider the blow-up along the sequence $X^k$ at scale $r_k=d(Y^k-X^k)$, i.e.
$$u^k(x,t)=\frac{1}{f(X^k)}\frac{u(x^k + r_k x, t^k+r_k^2 t)}{r_k^2} \to u^\infty(x,t)=P_{X_0}(x,t)$$
with 
$$P_{X_0}(x,t)=\frac12 {}^t x\cdot Q_{X_0}\cdot x + m_{X_0}t \quad \mbox{for}\quad (x,t)\in \R^n\times \R^-$$
and
$$Z^k=\left(\frac{y^k-x^k}{r_k},\frac{s^k-t^k}{r_k^2}\right) \to Z=(z,\tau),$$
with $Z\in \partial \left\{u^\infty>0\right\}$, $d(Z)=1$ and $\tau<0$. This gives a contradiction to the fact that $Q_{X_0}\ge 0$ and $m_{X_0}\le -\delta <0$.

\noindent {\bf Step 2: Conclusion}\\
Let a point $X_1\in \Gamma_\delta(n+0)$. Notice that (\ref{eq::el30}) implies that $\Gamma_\delta(n+0)$ is in particular a graph $t=h(x)$
in a small neighborhood of $X_1$.
By the classical Whitney's extension theorem (see \cite{Whi34}, or see Theorem 4 page 177 and paragraph 4.6 page 194 in \cite{Ste70}),
there exists a $C^2$ function (in the euclidean setting) $t=\tilde{h}(x)$ defined on $\R^n$, such that $\Gamma_\delta(n+0)$ is contained in the graph of $\tilde{h}$,
locally in a neighborhood of $X_1$.
\qeda

\begin{rem}
A finer result could be obtained in Proposition \ref{prop:graphn}, 
for the subset of singular points of $\Gamma_\delta(n+0)$ such that the matrix $Q_{X_0}$ 
has at least $j$ eigenvalues bounded from below by a positive constant 
for some $j\in \left\{0,...,n\right\}$.
\end{rem}

\noindent {\bf Proof of Theorem \ref{thm:mainsing2}}\\
First we observe that if $X_0\in \Gamma(k+1)$ (resp. $X_0\in \Gamma(n+0)$), then there are $\delta=\delta(X_0),r=r(X_0)>0$ 
such that $Q_r^-(X_0)\cap \Gamma_s\subset \Gamma_\delta(k+1)$  (resp. $Q_r^-(X_0)\cap \Gamma_s\subset \Gamma_\delta(n+0)$).

We deduce in particular that part 1 and part 2 of Theorem \ref{thm:mainsing2} 
follow respectively from Corollary \ref{cor:graph}, Proposition \ref{prop:graphn} and from the fact that $\Gamma_s$ is a closed set.
Part 3 of Theorem \ref{thm:mainsing2} now follows from the fact that every point in $\Gamma_s$
is either in $\Gamma(n+0)$ or in some $\Gamma(k+1)$ because we have
$$\Gamma_s = \Gamma(n+0)\cup \left(\bigcup_{k=0,...,n-1} \Gamma(k+1)\right)$$
and finally from the fact that for every $j\in \left\{0,...,n-1\right\}$, 
a  $(j+1)$-dimensional $C^1$-manifold is contained in a $((n-1)+1)$-dimensional $C^1$-manifold.
\qeda

\section{Pointwise decay estimates at singular points}\label{sec:decay}

In this section we obtain a decay rate of the distance between the solution 
and the unique blow-up at singular points. 
We first need to gather some estimates in order to obtain compactness.

\subsection{A global compactness result}

\begin{lem}{\bf (Global Cacciopoli-type inequality)}\label{lem:globalcacc}\\
Let $u$ be a solution of 
\begin{equation*}
\left\{\begin{array}{l}
Hu=h\chi_{\{u>0\}}\\
u\geq 0
\end{array}\textup{in $S_1$}\right.,
\end{equation*}
where $u$ has compact spatial support, $h\in L^2 (S_1)$ with $h(0)>0$ and $P/h(0)\in \ps$. 
Furthermore, put $w=u-P$. Then if $\eta\in C^\infty(I)$ with $I=[I_-,I^+]\subset (-1,0)$ there holds
\begin{align*}
\int_{I}\int_{\R^n}\left(|\nabla w|^2\eta+h(0)(\chi_{\{u>0\}}-\chi_{\{P>0\}})w\eta\right)\d\gamma+\frac12\int_{\R^n} w^2\eta \d\gamma^t\Big|_{t=I_-}^{t=I^+}\\
=\int_{I}\int_{\R^n}\left((h(0)-h)\chi_{\{u>0\}}w\eta+\frac12w^2\eta_t\right)\d\gamma.
\end{align*}
Here we recall that $\d\gamma^t=G(x,-t)\d t$.
\end{lem}

{\noindent \bf Proof of Lemma \ref{lem:globalcacc}}\\
We first start by noting that
$$(Hw)w=(h-h(0))\chi_{\{u>0\}}w+h(0)(\chi_{\{u>0\}}-\chi_{\{P>0\}})w.$$
By integrating by parts (which is justified as in Lemma \ref{lem:ipp}, recalling that $u$ has compact support, even if $w$ is not smooth but merely in $W_{2,\text{loc}}^{2,1}$) and using that $$\lap G(x,-t)=-\dd_t G(x,-t),$$ we obtain
\begin{align*}
\int_{\R^n}\int_{I}(Hw)w\eta \d\gamma & =\int_{\R^n}\int_{I}(-w_t+\lap
w)w\eta\d\gamma \\&=\int_{\R^n}\int_{I}\frac{\d}{\d t}(-\frac12
w^2\eta G(x,-t))\d x\d t +\int_{\R^n}\int_{I}\left(-|\nabla w|^2\eta+\frac12w^2\eta_t \right)\d\gamma .
\end{align*}
Putting these two equalities together we obtain
\begin{align*}
&\int_{\R^n}\int_{I}\left( |\nabla w|^2\eta+h(0)(\chi_{\{u>0\}}-\chi_{\{P>0\}})w\eta\right) \d\gamma+\frac12\int_{\R^n}
w^2\eta \d\gamma^t\Big|_{t=I_-}^{t=I^+}\\&=\int_{\R^n}\int_{I}\left((h(0)-h)\chi_{\{u>0\}}w\eta+\frac12w^2\eta_t\right)\d\gamma.
\end{align*}
\qeda

As a consequence we obtain the two following corollaries.
\begin{cor}{\bf ($L_t^\infty L_x^2$- and $L^2_tW^{1,2}_x$-estimates)}\label{cor:linftyw12}\\
Assume the hypotheses of Lemma \ref{lem:globalcacc}. Then for any $\delta\in (0,1)$ we have with $I_\delta=(-1+\delta,0)$
$$\int_{\R^n}\int_{I_\delta}\left(|\nabla w|^2+h(0)(\chi_{\{u>0\}}-\chi_{\{P>0\}})w\right)\d\gamma
+\sup_{t\in I_\delta}\int_{\R^n\times \{t\}}w^2\d\gamma^t\leq C(\delta,\sigma^h_2(1),\|w\|_{L^2_\gamma(S_1)}).$$ 
\end{cor}

{\noindent \bf ~Proof of Corollary \ref{cor:linftyw12}}\\
Take $t_0\in (-1+\delta,0)$ and use Lemma \ref{lem:globalcacc} with $0\leq \eta\leq 1$, $\eta(t_0)=1$, $\eta(-1+\delta/2)=0$ and $I=[-1+\delta/2,t_0]$. This implies
$$\int_{\R^n\times \{t_0\}}w^2\d\gamma^{t_0}\leq C(\delta,\sigma^h_2(1),\|w\|_{L^2_\gamma(S_1)}).$$
In order to estimate the remaining terms we choose $0\leq \eta\leq 1$, $\eta=1$ on $I_\delta^\e=[-1+\delta,-\e]$, $\eta(-1+\delta/2)=0$ and $I=[-1+\delta/2,-\e]$. Then Lemma \ref{lem:globalcacc} implies
$$\int_{\R^n}\int_{I_\delta^\e}\left(|\nabla w|^2+h(0)(\chi_{\{u>0\}}-\chi_{\{P>0\}})w\right)\d\gamma\leq C(\delta,\sigma^h_2(1),\|w\|_{L^2_\gamma(S_1)}).$$
Passing $\e\to 0$ yields the desired result.
\qeda

\begin{cor}{\bf ($L^1$-estimate)}\label{cor:l1}\\
Assume the hypotheses of Lemma \ref{lem:globalcacc}. Then for any $\delta\in (0,1)$ we have with $I_\delta=(-1+\delta,0)$
$$\int_{\R^n}\int_{I_\delta}h(0)|\chi_{\{u>0\}}-\chi_{\{P>0\}}|\d\gamma\leq C(\delta,\sigma^h_2(1),\|w\|_{L^2_\gamma(S_1)}).$$ 
\end{cor}

{\noindent \bf Proof of Corollary \ref{cor:l1}}\\
We have
$$(Hw)\sgn w=(h-h(0))\chi_{\{u>0\}}\sgn w +h(0)|\chi_{\{u>0\}}-\chi_{\{P>0\}}|.$$
By partial integration (again justified as in Lemma \ref{lem:ipp}) and using that $(\sgn (t))'\geq 0$, $\lap G(x,-t)=-\dd_t G(x,-t)$ and a regularization of the sign function, we deduce
$$\int_{\R^n}\int_{I_\delta^\e} (Hw)\sgn w\d\gamma\leq \int_{\R^n}\int_{I_\delta^\e}\dd_t(-|w|G(x,-t)) \d x \d t,$$
where $I_\delta^\e = [-1+\delta, -\e]$. Putting this together we arrive at
\begin{align*}
&\int_{\R^n}\int_{I_\delta^\e} h(0)|\chi_{\{u>0\}}-\chi_{\{P>0\}}|\d\gamma\\&\leq
\int_{\R^n}\int_{I_\delta^\e}\dd_t(-|w|G(x,-t))\d x \d t+\int_{\R^n}\int_{I_\delta^\e}(h(0)-h)\chi_{\{u>0\}}\sgn
w\d\gamma\\&\leq C(\delta,\sigma^ h_2(1),\|w\|_{L^2_\gamma(S_1)}),
\end{align*}
by Corollary \ref{cor:linftyw12} and H\"older's inequality. To conclude the proof we let $\e\to 0$.
\qeda

Below we state Simon's compactness theorem (see Corollary 6 in \cite{Sim87}).

\begin{thm}\label{thm:aubin}{ \bf (Compactness in Banach spaces)}\\
Let $X_0\subset X\subset X_1$ be Banach spaces such that $X_0$ is
compactly embedded in $X$ and $X$ is continuously embedded in
$X_1$. Moreover, assume that $u_k$ is a sequence of functions such that for some $q>1$ 
$$||u_k||_{L^q(I;X)}+||u_k||_{L^1(I;X_0)}+||\dd_t u_k||_{L^1(T;X_1)}\leq C,$$
where $I\subset \R $ is a compact interval. Then there is a subsequence
$u_{k_j}$ that converges in $L^p(I;X)$ for all $1\leq p<q$.
\end{thm}

In the lemma that follows, we show that the spaces we wish to use
are good for Simon's theorem.

\begin{lem}{\bf (Banach space inclusions)}\label{lem:inclusions}\\
\begin{enumerate}
\item Let $Y$ be the closure of $C_0^\infty(\R^n)$ with respect to the norm
$$\|\phi\|_{Y}= \Big\|\frac{\phi}{G(x,1)}\Big\|_{L^\infty(\R^n)}+\Big\|\frac{\nabla \phi}{\sqrt{G(x,1)}}\Big\|_{L^2_{(\R^n)}}.$$
Then for any $u=\nabla v+w$ with $v\in L^2_{\gamma^{-1}}(\R^n)$ and $w\in L^1_{\gamma^{-1}}(\R^n)$ we have
\begin{equation}\label{eq:y1}
\|u\|_{Y'}\leq \|v\|_{L^2_{\gamma^{-1}}(\R^n)}+\|w\|_{L^1_{\gamma^{-1}}(\R^n)},
\end{equation}
and in particular
\begin{equation}\label{eq:y2}
\|u\|_{Y'}\leq C\|u\|_{L^2_{\gamma^{-1}}(\R^n)},
\end{equation}
for any $u\in L^2_{\gamma^{-1}}(\R^n)$. Here we recall that
$$
\|u\|^p_{L^p_{\gamma^{-1}}(\R^n)}=\int_{\R^n} |u|^p d\gamma^{-1} = \int_{\R^n} |u|^p G(x,1) d x.
$$
\item The embedding $H^1_{\gamma^{-1}}(\R^n)\to L^2_{\gamma^{-1}}(\R^n)$ is compact. Here we use the notation
$$
\|u\|_{H^1_{\gamma^{-1}}(\R^n)}^2=\int_{\R^n} \left(u^2+|\nabla u|^2\right) d\gamma^{-1}.
$$
\end{enumerate}
\end{lem}

{\noindent \bf Proof of Lemma \ref{lem:inclusions}}\\
To prove the first part take $\phi\in C_0^\infty(\R^n)$. Then
\begin{align*}
|\<\phi,u\>|& =|\<-\nabla \phi,v\>+\<\phi,w\>|\leq \int_{\R^n}| v\nabla \phi|+|\phi w| \d x\\
&\leq \int_{\R^n}|\nabla \phi|G(x,1)^{-\frac12}|
v|G(x,1)^{\frac12}+|\phi G(x,1)^{-1} w G(x,1)|\d x\\
&\leq \|\phi\|_Y\left(\|v\|_{L^2_{\gamma^{-1}}(\R^n)}+\|w\|_{L^1_{\gamma^{-1}}(\R^n)}\right),
\end{align*}
by H\"older's inequality.

Now, for any $u\in L^2_{\gamma^{-1}}(\R^n)$ we have from H\"older's inequality
$$
\|u\| _{L^1_{\gamma^{-1}}(\R^n)}\leq \|u\| _{L^2_{\gamma^{-1}}(\R^n)}\big \|\sqrt{G(x,1)}\big \|_{L^2(\R^n)}\leq C\|u\| _{L^2_{\gamma^{-1}}(\R^n)}.
$$
This combined with \eqref{eq:y1} implies the estimate \eqref{eq:y2}.

The second part follows from Theorem 3.1 in \cite{Hoo81}.
\qeda

\begin{lem}{\bf (Global compactness)}\label{lem:compactglobal}\\
Assume we have two sequences of functions $u_k$ and $P_k$ such that
$$\left\{\begin{array}{l}
\begin{array}{l} 
H u_k=f_k\chi_{\{u_k>0\}}\\ 
u_k\geq 0\\
\end{array} \Big|\textup{ in $S_2$},  \\
\text{the support of $u_k$ is contained in $B_{R_k}\times (-4,0)$ for some $R_k>0$},\\
\displaystyle\frac{P_k}{f_k(0)}\in \ps.
\end{array}\right. $$
Assume further that with $w_k=u_k-P_k$ there holds
\begin{equation}\label{eq:cest}
||w_k||_{L_\gamma^2(S_2)}+\sigma_2^{f_k}(2)\leq C.
\end{equation}
Then, for each $\delta\in (0,1)$ there is a subsequence of $w_k$ converging in $L_\gamma^2(\R^n\times (-1,-\delta))$.
\end{lem}

{\bf \noindent Proof of Lemma \ref{lem:compactglobal}}\\
Using \eqref{eq:cest},  we can apply Corollary \ref{cor:linftyw12} and Corollary \ref{cor:l1} to deduce 
\begin{equation}\label{eq:wkest}
\|w_k\|_{L^\infty_t(L^2_{\gamma^t,x})(S_1)},\|\nabla w_k\|_{L^2_\gamma(S_1)},\|f_k(0)|\chi_{\{u_k>0\}}-\chi_{\{P_k>0\}}|\|_{L^1_\gamma (S_1)}\leq C,
\end{equation}
where we use the notation
$$
\|u\|^q_{L^q_t(L^p_{\gamma^t,x})(S_1)}=\int_{(-1,0)}\|u\|^q_{L^p_{\gamma^t}(\R^n)}d t=\int_{(-1,0)}\left(\int_{\R^n} |u(x,t)|^p G(x,-t) dx \right)^\frac{q}{p} dt,
$$
with the standard convention for the case $q=\infty$.
Now define $\hat{w}_k(x,t)=w_{k}(x\sqrt{-t},t)$. Then from \eqref{eq:wkest} and a change of variables
$$
\int_{\R^n}\hat w_k^2(x,t) G(x,1)\d x=\int_{\R^n} w_k^2 \d\gamma^t\leq C
$$
and
$$
\int_{S_1} |\nabla \hat w_k(x,t)|^2 G(x,1)\d x\d t\leq \int_{S_1} |\nabla w_k(x,t)|^2|t|\d\gamma\leq C.
$$
Thus, 
\begin{equation}\label{eq:simon1}
\|\hat w_k\|_{L^\infty_t(L^2_{\gamma^{-1},x})(S_1)},\|\nabla\hat  w_k\|_{L^2_t(L^2_{\gamma^{-1},x})(S_1)}\leq C.
\end{equation}
In addition
\begin{align*}
\dd_t \hat w_k(x,t)&=(\dd_t w_k)(x\sqrt{-t},t)-\frac{x}{2\sqrt{-t}}\cdot \nabla w_k (x\sqrt{-t},t)\\
&= \left(\lap w_k+f_k(0)(\chi_{\{u_k>0\}}-\chi_{\{P_k>0\}})+(f_k-f_k(0))\chi_{\{u_k>0\}})\right)(x\sqrt{-t},t)\\&-\frac{x}{2\sqrt{-t}}\cdot \nabla w_k (x\sqrt{-t},t)\\
&=a_k+b_k+c_k+d_k.
\end{align*}
For the first term we have
$$
a_k = (\lap w_k)(x\sqrt{-t},t)=\nabla \cdot \left(\frac{1}{\sqrt{-t}}(\nabla w_k)(x\sqrt{-t},t)\right),
$$
where
$$
\int_{\R^n\times (-1,-\delta)}\frac{1}{(-t)}|(\nabla w_k)(x\sqrt{-t},t)|^2 G(x,1)\d x\d t\leq \frac{1}{\delta}\|\nabla w_k\|_{L^2_\gamma(S_1)}^2\leq C_\delta, 
$$
for any $\delta\in (0,1)$, by \eqref{eq:wkest}. From \eqref{eq:wkest} we can also conclude
$$
\int_{S_1} |b_k|G(x,1)\d x \d t\leq \int_{S_1} f_k(0)|\chi_{\{u_k>0\}}-\chi_{\{P_k>0\}}|\d\gamma \leq C,
$$
whereas the estimate
$$
\int_{S_1} |c_k|G(x,1)\d x\d t\leq \int_{S_1}|(f_k-f_k(0))\chi_{\{u_k>0\}}|\d\gamma \leq C\sigma_2^{f_k}(2)\leq C,
$$
follows from \eqref{eq:cest}. Moreover, 
$$
\int_{\R^n\times (-1,-\delta)} |d_k|G(x,1)\d x\d t\leq \Big\|\frac{x}{2t}\Big\|_{L^2_\gamma(\R^n\times (-1,-\delta))}\|\nabla w_k\|_{L^2_\gamma(S_1)}\leq C_\delta,
$$
for any $\delta\in (0,1)$, by H\"older's inequality and \eqref{eq:wkest}. We can thus draw the conclusion that
$$
\dd_t \hat w_k = a_k+(b_k+c_k+d_k)=\nabla e_k+(b_k+c_k+d_k)
$$
with
$$
\|e_k\|_{L^2_t(L^2_{\gamma^{-1},x})(\R^n\times (-1,-\delta))}+ \|b_k+c_k+d_k\|_{L^1_{\gamma^{-1}}(\R^n\times (-1,-\delta))}\leq C_\delta,
$$
for any $\delta\in  (0,1)$. Hence, from the first part of Lemma \ref{lem:inclusions}
\begin{equation}\label{eq:simon2}
\|\dd_t \hat{w_k}\|_{L^1_t(Y')(\R^n\times (-1,-\delta))}\leq C_\delta,
\end{equation}
for any $\delta\in (0,1)$.
Applying Simon's Theorem (Theorem \ref{thm:aubin}) with $I=(-1,-\delta)$, $X_0=H^1_{\gamma^{-1}} (\R^n)$, $X=L^2_{\gamma^{-1}}
(\R^n)$ and $X_1=Y'$, we obtain from \eqref{eq:simon1} and \eqref{eq:simon2} that there is a subsequence of $\hat{w}_k$
converging in $L^p_t(L^2_{\gamma^{-1},x})(\R^n\times (-1,-\delta))$ for all $1\leq p<\infty$, 
and in particular it converges in $L^2_{\gamma^{-1}}(\R^n\times (-1,-\delta))$. 
Transforming back to $w_k$ we see that there is a subsequence of $w_k$ converging in $L^2_\gamma(\R^n\times (-1,-\delta))$.
\qeda




\subsection{Liouville-type results}

In this section, we present some technical results needed for the proofs
given later. The first result, is the parabolic counterpart to Lemma 7.6 in \cite{Mon09}.

\begin{lem}{\bf (Liouville-type result (I))}\label{lem:homocaloric}\\
Let $u\in L^2_\gammab(S_1)$ satisfy
\begin{enumerate}
\item $u$ is parabolically homogeneous of degree two.
\item $Hu=0$ in $(\R^n\times (-\infty,0))\cap \{x_1\neq 0\}$, 
\item $Hu\leq 0$ in $\R^n\times (-\infty,0)$,
\end{enumerate}
Then $Hu=0$ in $\R^n\times (-\infty,0)$ and $u$ must be a caloric polynomial of degree two.
\end{lem}

\noindent {\bf Proof of Lemma \ref{lem:homocaloric}}\\
Let $P=P(x)$ be a caloric function, i.e.
\begin{equation}\label{eq::el53}
H P=0
\end{equation}
such that $P$ is parabolically homogeneous of degree two. 
Due to the homogeneity of $P$ and $u$ we have
\begin{equation}\label{eq::el50}
x\cdot \nabla P=2P- 2tP_t,
\end{equation}
\begin{equation}\label{eq::el51}
x\cdot \nabla u=2u-2tu_t
\end{equation}
Take $\phi(x)\ge 0$ be a smooth function such that
$$
\phi =1\textup{ on $B_1$},\quad \phi =0 \textup{ outside $B_2$},
$$
and assume moreover that $\phi$ is radial and nonincreasing in $|x|$.
Let $\phi_R(x)=\phi(x/R)$. Take $\psi(t)\ge 0$ be a smooth function with support contained $I\subset (-1,0)$, and set
$$\zeta_R(x,t)=\phi_R(x) \psi(t).$$
Now let $\check{G}(x,t)=G(x,-t)$ and recall that 
\begin{equation}\label{eq::el52}
\nabla \check{G}(x,t)=\frac{x}{2t}\check{G}(x,t).
\end{equation}
We integrate by parts, using \eqref{eq::el53}, \eqref{eq::el50}, \eqref{eq::el51} and \eqref{eq::el52}, to obtain
\begin{align*}
\langle \Delta u, \zeta_R \check{G}P\rangle &
= - \langle \nabla u, \nabla (\zeta_R \check{G}P)\rangle \\
& = - \langle \nabla u, (\nabla \zeta_R) \check{G}P + \zeta_R (\nabla \check{G})P + \zeta_R \check{G} (\nabla P)\rangle \\
&= - \langle \nabla u, (\nabla \zeta_R) \check{G}P + \zeta_R (\nabla \check{G})P\rangle  
+ \langle u, (\nabla \zeta_R \cdot \nabla P) \check{G}  + \zeta_R (\nabla\check{G} \cdot \nabla P) + \zeta_R \check{G} (\Delta P) \rangle \\
& = - \langle \nabla u, (\nabla \zeta_R) \check{G}P + \zeta_R \left(\frac{x}{2t} \check{G}\right)P\rangle  \\
&+ \langle u, (\nabla \zeta_R \cdot \nabla P) \check{G}  + \zeta_R \left(\check{G}\ \frac{x}{2t} \cdot \nabla P\right)
+ \zeta_R \check{G}P_t\rangle \\
&= - \langle \nabla u, (\nabla \zeta_R) \check{G}P\rangle  -\langle \frac{x}{2t}\cdot \nabla u, \zeta_R \check{G} P\rangle  
+ \langle u, (\nabla \zeta_R \cdot \nabla P) \check{G}\rangle   + \langle u,\zeta_R \check{G}\ \frac{P}{t} \rangle \\
&= - \langle \nabla u, (\nabla \zeta_R) \check{G}P\rangle  + \langle u_t, \zeta_R \check{G} P\rangle  
+ \langle u, (\nabla \zeta_R \cdot \nabla P) \check{G}\rangle.
\end{align*}
In other words, we have, using \eqref{eq::el52} again
\begin{align*}
\langle Hu, \zeta_R \check{G}P\rangle & =\langle u, (\nabla \zeta_R \cdot \nabla P) \check{G}\rangle - \langle \nabla u, (\nabla \zeta_R) \check{G}P\rangle \\
&=\langle u, (\nabla \zeta_R \cdot \nabla P) \check{G}\rangle 
+ \langle u, (\Delta \zeta_R) \check{G}P + (\nabla \zeta_R \cdot\nabla \check{G})P + (\nabla \zeta_R\cdot \nabla P) \check{G}\rangle \\
&=\langle u, (\Delta \zeta_R) \check{G}P + (\nabla \zeta_R \cdot\nabla \check{G})P + 2(\nabla \zeta_R\cdot \nabla P) \check{G}\rangle \\
&=\langle u, K_R\check{G}\rangle,
\end{align*}
where we introduced
$$
 K_R:=\left((\Delta \phi_R)P + (\nabla \phi_R \cdot \frac{x}{2t})P + 2\nabla \phi_R\cdot \nabla P\right)\psi.
 $$
Using Young's inequality and that $\supp \nabla \phi_R\subset B_{2R}\backslash B_R$ we can then deduce
$$\langle -Hu, \zeta_R \check{G}P\rangle  = -\int_{\R^n\times (-\infty,0)} 
u K_R d\gamma \le \int_{(B_{2R}\backslash B_R)\times I} (u^2 + K_R^2) d\gamma \longrightarrow  0$$
as $R\to \infty$, since $u\in L^2_\gammab(S_1)$ and since the Gaussian is integrable against any polynomial. 

With the choice
$$P(x,t)=-n x_1^2+\left(\sum_{i=2}^{n}x_i^2 -2t\right),$$
we have 
\begin{equation}\label{eq::el56}
P> 0\quad \mbox{on}\quad \left\{x_1=0\right\}\subset \R^n\times (-\infty,0). 
\end{equation}
By assumption (3), $\mu=-Hu$ is a non-negative measure on $\R^n\times (-\infty,0)$
and with its support contained in $\left\{x_1=0\right\}$. Thus, 
$\langle \mu,\zeta_R \check{G}P\rangle$ is non-negative, monotonically nonincreasing in $R$ and tends to zero as $R\to \infty$. 
The only possibility is that $\langle \mu,\zeta_R \check{G}P\rangle$ is zero for all $R>0$.
From (\ref{eq::el56}), we deduce that $\mu=0$, i.e., $u$ is caloric in $\R^n\times (-1,0)$ and then also in $\R^n\times (-\infty,0)$ by the homogeneity.

Since $u$ is caloric in for instance $Q_2^-$, we know that $u$ is smooth in $Q_{1}^-$. In particular, $u$ is in $L^\infty(Q_{1}^-)$ (but we do neither need nor have an explicit control of the bound). Interior estimates (Theorem \ref{thm:interior}), the parabolic Sobolev embedding (Theorem \ref{thm:sobolev}) and the homogeneity of $u$ imply that for any $i,j,k=1,\ldots,n$ and any $R>1$
$$
\sup_{Q_R^-}\left(|\dd_{i,j,k} u|+|\nabla u_t|+|u_{tt}|\right)\leq \frac{C}{R^3}\|u\|_{L^\infty(Q_{2R}^-)}=\frac{4CR^2}{R^3}\|u\|_{L^\infty(Q_1^-)}\to 0, 
$$
as we let $R\to\infty$. Thus, $\dd_{i,j,k} u = |\nabla u_t|=u_{tt}=0$ and $u$ must be a parabolic polynomial of degree two.
\qeda\\

Next we present a Liouville-type result which might appear quite
strange at first glance, but it turns out that it is exactly the result
we will need later on.

\begin{prop}{\bf (Liouville-type result (II))}\label{prop:liouville}\\
Assume
\begin{enumerate}
\item $Hv=0$ in $\R^n\times \R^-$, 
\item $v$ is parabolically homogeneous of degree 2, 
\item $v\geq 0$ on $\{x_i=0, i=1, \ldots, k, t=0\}$ for some integer $k\in [0,n]$,
\item 
$$\int_{S_1}vx_ix_j\d\gammab=0$$
for $i=1,\ldots,k, j=1,\ldots,n$ and $i\neq j$.
\item 
$$\int_{S_1}v q\d\gammab\leq 0$$
for all functions $p$ of the form
$$p=\sum_{i=1}^n\gamma_ix_i^2+\alpha t$$
where 
$$
\begin{array}{lr}
\displaystyle\sum_{i=1}^n 2\gamma_i=\alpha,\\
\gamma_i\in \R\textup{ for }i\leq k,\\
\gamma_i\in \R^+\textup{ for }i>k,
\end{array}
$$
and either\\
({\bf Case 1}):
$$
\alpha\geq 0,
$$
or ({\bf Case 2}):
$$
\begin{array}{lr}
\alpha\leq 0,\\
v\geq 0\text{ on }\{x_i=0,i=1,\ldots,k,t\leq 0\}.
\end{array}
$$
\end{enumerate}
Then $v=0$ in $\R^n\times \R^-$.
\end{prop}

{\noindent \bf Proof of Proposition \ref{prop:liouville}}\\
Since $v$ is caloric and parabolically homogeneous of degree two, we can as in the proof of Lemma \ref{lem:homocaloric} conclude that $v$ is a caloric polynomial homogeneous of degree two. We thus can write
$$v=\frac12 \sum_{i,j\leq n}q_{ij}x_ix_j+mt$$
with $(Q)=(q_{ij})$ a symmetric matrix and $\tr Q=m$. 
Property (3) together with property (4) implies that
$$Q=\left(\begin{array}{lr}
A& 0\\
0& B
\end{array}\right),$$
with $A$ a diagonal $k\times k$-matrix and $B$ a positive semi-definite $(n-k)\times (n-k)$-matrix.
\noindent {\bf Step 1: Proving \eqref{eq:gamma}.}\\ \noindent 
We use property (5) with
$$p=\sum_{i=1}^k\gamma_Ax_i^2+\sum_{i=k+1}^n\gamma_Bx_i^2+\alpha t$$ 
and 
\begin{equation}\label{eq:gammarel}
2k\gamma_A+2(n-k)\gamma_B=\alpha,\quad \gamma_A\in \R, \gamma_B\geq 0.
\end{equation}
Then (5) reads
\begin{equation}\label{eq:vqleq}
0\geq \int_{S_1}vp\d\gammab=\sum_{i=1}^6I_i,
\end{equation}
where
\begin{align*}
I_1&=\frac12\int_{S_1}
\sum_{i,j\leq n}q_{ij}x_jx_i\sum_{l=1}^k\gamma_Ax_l^2\d\gammab\\
&=\frac{\gamma_A}{2}\int_{S_1}\sum_{i\leq n,l\leq k}q_{ii}x_i^2x_l^2\d\gammab=\frac{\gamma_A}{2}\int_{S_1}\sum_{l=1}^k\left(q_{ll}x_l^4+\sum_{i\neq
    l}q_{ii}x_i^2x_l^2\right)\d\gammab\\
&=\frac{\gamma_A}{2}\left(\sum_{l=1}^k(q_{ll}a+\sum_{i\neq
    l}q_{ii}b)\right)
=\frac{\gamma_A}{2}a\tr
A+\frac{\gamma_A}{2}b\sum_{l=1}^k\left(\sum_{i=1}^nq_{ii}-q_{ll}\right)\\
&=\frac{\gamma_A}{2}\left(a\tr A+bk\tr Q-b\tr A\right)=\frac{\gamma_A}{2}\left((a-b)\tr
  A+kb\tr Q\right),
\end{align*}
\begin{align*}
I_2&=\frac12\int_{S_1}
\sum_{i,j\leq n}q_{ij}x_jx_i\sum_{l=k+1}^n\gamma_Bx_l^2\d\gammab=
\frac{\gamma_B}{2}\left((a-b)\tr
  B+(n-k)b\tr Q\right),\\
I_3&=\int_{S_1}\frac12\sum_{i,j\leq  n} q_{ij}x_i x_j\alpha t\d\gammab=\frac{\alpha c}{2}\sum_{i=1}^n
q_{ii}=\frac{\alpha c}{2}\tr Q ,\\
I_4&=\int_{S_1} \sum_{j=1}^k\gamma_A x_j^2mt\d\gammab=m\gamma_Akc,\\
I_5&=\int_{S_1} \sum_{j=k+1}^n\gamma_B x_j^2mt\d\gammab=m(n-k)\gamma_Bc,\\
I_6&=\int_{S_1} m\alpha t^2\d\gammab=m\alpha d,\\
\end{align*}
where we have used the notation
$$a=\int_{S_1}\frac{x_1^4}{-t}\d \gamma,\quad 
b=\int_{S_1}\frac{x_1^2x_2^2}{-t}\d\gamma ,\quad c=\int_{S_1}\frac{t x_1^2}{-t}\d \gamma\quad\textup{and}\quad d=\int_{S_1}\frac{t^2}{-t}\d\gamma.$$
Using the relation between $\gamma_A,\gamma_B$ and $\alpha$, and that $\tr Q=m$, we obtain
\begin{align}
\nonumber \sum_{i=1}^6 I_i&=
\frac{\gamma_A}{2}\left((a-b)\tr
  A+kb m\right)+\frac{\gamma_B}{2}\left((a-b)\tr
  B+(n-k)b m\right)\\\nonumber &+\frac{\alpha}{2}mc+m\gamma_Akc+m(n-k)\gamma_Bc+m\alpha
d\\\label{eq:sumi}
&=\frac{\gamma_A}{2}(a-b)\tr A+\frac{\gamma_B}{2}(a-b)\tr B+\frac{bm\alpha}{4}+\frac{mc\alpha}{2}+\frac{\alpha}{2}mc +m\alpha
d\\
\nonumber &=\frac{a-b}{2}(\gamma_A\tr A+\gamma_B \tr
B)+m\alpha(\frac{b}{4}+c+d).
\end{align}
Up to the multiplicative constant
$$\int_0^1(-t)\int_{\R^n}e^{-\frac{y^2}{4}}\d y \d t>0,$$
we find $a=12,b=4,c=-2$ and $d=1$, and thus
$$\frac{b}{4}+c+d=0.$$
Then \eqref{eq:vqleq} and \eqref{eq:sumi} imply
$$\frac{a-b}{2}(\gamma_A\tr A+\gamma_B \tr B)\leq 0,$$
which is equivalent to
\begin{equation}\label{eq:gamma}
\gamma_A\tr A+\gamma_B \tr B=(\gamma_B-\gamma_A)\tr B+\gamma_A m\leq 0.
\end{equation}
\noindent {\bf Step 2: Conclusion.}\\ \noindent 
Now we divide the proof into different cases. First we choose $\alpha=1$.


{\bf Case a: ($k\neq 0, k\neq n$)} \\
We can choose $\gamma_A$ 
and $\gamma_B$ to be whatever we need as long as \eqref{eq:gammarel} holds. 
In particular we can take $\gamma_B>\gamma_A$ and $\gamma_A$ having the
same sign as $m$. Then \eqref{eq:gamma} implies $\tr B=0$ and $m=0$, since $\tr B\geq 0$.

{\bf Case b: ($k=0$)} \\
Taking $\gamma_B>0$ in \eqref{eq:gamma} yields $\tr B\leq 0$ and thus $0=\tr B=\tr Q =m$ in this case.

{\bf Case c: ($k=n$)} \\ We can choose $\gamma_A>0$ and then, since we have no $B$-terms, \eqref{eq:gamma} simply implies $\tr A=m\leq 0$.





In order to conclude that $m=0$ also in the case $k=n$, we now apply the second alternative (Case 2), by choosing $\gamma_A<0$ and $\alpha=\alpha(\gamma_A)<0$. Plugged into \eqref{eq:gamma} this implies $m=\tr A = 0$.

Hence, all the cases above, we can deduce that $\tr B=0$ and $\tr A=m= 0$ . Since $B$ is positive semi-definite $B=0$. With
$$p=\sum_{i=1}^kq_{ii}x_i^2,$$
we obtain from (5) that
$$\int_{S_1}\sum_{i,j=1}^kq_{ii}q_{jj}x_i^2x_j^2\d\gammab=b(\tr A)^2+(a-b)\sum_{i=1}^kq_{ii}^2\leq 0,$$
which forces $A=0$.
\qeda\\

The lemma below is a quite standard result on minimizing polynomials,
partially already proved in \cite{Mon09}.

\begin{lem}{\bf (Estimates of $\ns$ in smaller and bigger strips)}\label{lem:dyadic}\\
Suppose
$$\ns(u,1)=\|u-P_1\|_{L^2_\gammab(S_1)},$$
where $P_1\in \mathcal{P}_\textup{sing}$. Then for any $\rho\geq 1$
$$\left(\frac{1}{\rho^4}\int_{S_\rho}|u-P_1|^2\d\gammab\right)^\frac12\leq C\int_1^{2\rho}\frac{\ms(u,s)}{s}\d s$$
and for $\rho\leq 1$
$$\left(\frac{1}{\rho^4}\int_{S_\rho}|u-P_1|^2\d\gammab\right)^\frac12\leq C\int^3_{2\rho}\frac{\ms(u,s)}{s}\d s.$$
\end{lem}

{\noindent\bf Proof of Lemma \ref{lem:dyadic}}\\
For the case  $\rho\geq 1$, the proof is exactly the same as in the proof of Lemma 2.9 in \cite{Mon09}. 
We give the proof for the case $\rho \leq 1$, which also is almost identical with that proof. 
Throughout the proof we write $M$ instead of $\ms$, and we denote by $P_r$ the element in $\ps$ corresponding to $\ns (u,r)$. We have
\begin{align}
\nonumber &\left(\frac{1}{r^4}\int_{S_r}|P_{2r}-P_r|^2\d
  \gammab\right)^\frac12\\ &\leq
\left(\frac{1}{r^4}\int_{S_r}|u-P_r|^2\d
  \gammab\right)^\frac12+\left(\frac{1}{r^4}\int_{S_r}|u-P_{2r}|^2\d
  \gammab\right)^\frac12 \label{eq:dyadicest1}\\
&\leq M(r)+2^2M(2r)\leq CM(2r).\nonumber
\end{align}
By the same reasoning, we can for fixed $r\in (1/2,1]$ obtain
\begin{equation}
\left(\frac{1}{r^4}\int_{S_r}|P_{r}-P_1|^2\d \gammab\right)^\frac12\leq C
M(1).\label{eq:dyadicest2}
\end{equation}
Now we remark that for any $\rho>0$ there holds
\begin{align}
\nonumber &\left(\frac{1}{\rho^4}\int_{S_\rho}|u-P_1|^2\d
  \gammab\right)^\frac12\\
&\leq\left(\frac{1}{\rho^4}\int_{S_\rho}|u-P_\rho|^2\d
  \gammab\right)^\frac12+\left(\frac{1}{\rho^4}\int_{S_\rho}|P_\rho-P_1|^2\d
  \gammab\right)^\frac12\label{eq:rhoest} \\
&\leq M(\rho)+\left(\frac{1}{\rho^4}\int_{S_\rho}|P_\rho-P_1|^2\d
  \gammab\right)^\frac12\nonumber .
\end{align}
Therefore, if we choose $\rho=2^{-k}r$  and use \eqref{eq:dyadicest1}, \eqref{eq:dyadicest2}, \eqref{eq:rhoest} and the homogeneity of the $P_r$s, we deduce
\begin{align}&\nonumber \left(\frac{1}{\rho^4}\int_{S_\rho}|P_\rho-P_1|^2\d
  \gammab\right)^\frac12 \\&\leq \left(\frac{1}{\rho^4}\int_{S_\rho}|P_r-P_1|^2\d
  \gammab\right)^\frac12+\sum_{j=0}^{k-1}\left(\frac{1}{\rho^4}\int_{S_\rho}|P_{2^{-j}r}-P_{2^{-j-1}r}|^2\d
  \gammab\right)^\frac12\label{eq:pdyad}\\\nonumber 
&\leq C M(1)+C\sum_{j=0}^{k-1} M(2^{-j}r).
\end{align}
Thus, from \eqref{eq:rhoest} and \eqref{eq:pdyad} it follows
\begin{align*}
\left(\frac{1}{\rho^4}\int_{S_\rho}|u-P_1|^2\d
  \gammab\right)^\frac12&\leq C\sum_{j=0}^{k-1} M(2^{-j}r)+C M(1)+M(\rho)\\
&\leq C\sum_{j=0}^kM(2^{-j}2r)\\&\leq C\int_{2\rho}^3\frac{M(u,r)}{r}\d r.
\end{align*}
\qeda\\

\subsection{Decay estimate}

Now we can prove the following result which is the key ingredient in the proof of Theorem \ref{thm:mainsing}.

\begin{prop}\label{prop:decaysing}{\bf (Decay estimate at singular points)}\\
Let $u$ be a solution of \eqref{eq:obstacle} with $p>(n+2)/2$, $p\geq 2$ and $\tilde\sigma_p(1)$ finite. Then there are constants $M_0,C_0>0$, $r_0,\lambda,\mu\in (0,1)$ such that for all $r<r_0$ 
$$M_{\textup{sing}}(u,r)\leq M_0\quad \Rightarrow \left\{\begin{array}{lr} M_{\textup{sing}}(u,\lambda r)<\mu \ms(u,r)\\
\textup{or}\\
\ms(u,r)<C_0\Sigma^f(r)
\end{array}\right.$$
where
$$\Sigma^f(\tau)=\sigma^f_p(\tau)+\int_0^\tau \frac{\sigma^f_p(r)}{r}.$$
\end{prop}

{\noindent \bf Proof of Proposition \ref{prop:decaysing}}\\
The proof is by contradiction and it is divided into several steps. 
Throughout the proof we will write $M$ instead of $M_\textup{sing}$.

\noindent {\bf Step 1: Initialization of the blow-up procedure and gathering of estimates}.\\
If the assertion is not true  then there are
$$C_k\to\infty, M_k,r_k,\lambda_k\to 0\text{ and }\mu_k\to 1$$
such that the statement is false for a sequence of solutions $u_k$ to $Hu_k=f_k\chi_{\{u_k>0\}}$. In other words, $M(u_k,r_k)\leq M_k$ but still 
\begin{equation}\label{eq:muk}
M(u_k,\lambda_kr_k)\geq \mu_k
M(u_k,r_k)
\end{equation}
and
\begin{equation}\label{eq:Msigma}
M(u_k,r_k)\geq C_k\Sigma_{f_k}(r_k),
\end{equation}
where we hereafter in the proof will denote $\Sigma^{f_k}$ simply by $\Sigma$. We can also assume that for some $0<\rho_k\leq \lambda_kr_k$ there holds
\begin{equation}\label{eq:epsk}
\frac{M(u_k,\lambda_kr_k)}{1+1/k}\leq N(u_k,\rho_k)=\e_k\to 0.
\end{equation}
Let
$$v_k(x,t)=\frac{u_k(\rho_kx,\rho_k^2t)}{\rho_k^2}$$
and
$$w_k(x,t)=\frac{u_k(\rho_kx,\rho_k^2t)-P_k(\rho_kx,\rho_k^2t)}{\e_k\rho_k^2},$$
where $P_k\in \ps $ is a function that realizes the infimum of
$N$ at the level $\rho_k$. Now we have 
\begin{equation}\label{eq:T1est}
  \inf_{P\in \ps}\left(\int_{S_1}\left|w_k-\frac{P-P_k}{\e_k}\right|^2\d\gammab\right)^\frac12=1
\end{equation}
with equality for $P=P_k$. Moreover, for $s\rho_k\leq r_k$, \eqref{eq:muk} and \eqref{eq:epsk} implies
\begin{equation}\label{eq:lessthanone}
\inf_{P\in P_\textup{sing}}\left(\frac{1}{s^4}\int_{S_{s}}\left|w_k-\frac{P-P_k}{\e_k}\right|^2\d\gammab\right)^\frac12\leq \frac{1+1/k}{\mu_k}\to 1.
\end{equation}
Furthermore,
$$N(v_k,1)=\left(\int_{S_1}|v_k-P_k|^2\d\gammab\right)^\frac12=\e_k,$$
so by applying Lemma \ref{lem:dyadic} together with \eqref{eq:muk}, \eqref{eq:Msigma}, \eqref{eq:epsk} and using that $\rho_k\leq \lambda_k r_k$ with  $\lambda_k\to 0$, we obtain for $s\leq 1$ and $k$ large enough
\begin{align}\label{eq:subquadr}
&\nonumber\left(
  \frac{1}{s^4}\int_{S_s}|w_k|^2\d\gammab\right)^\frac12=\frac{1}{\e_k}\left(
  \frac{1}{s^4}\int_{S_s}|v_k-P_k|^2\d\gammab\right)^\frac12\\
&\leq
\frac{C}{\e_k}\int_{2s}^3\frac{M(v_k,\tau)}{\tau}\d\tau=\frac{C}{\e_k}\int_{2s}^3\frac{M(u_k,\rho_k\tau)}{\tau}\d\tau\\
\nonumber
&\leq \frac{C}{\e_k}\int_{2s}^3\frac{M(u_k,r_k)}{\tau}\d\tau \leq
\frac{C}{\e_k}\frac{1+1/k}{\mu_k}N(u_k,\rho_k)(|\ln 2s|+\ln 3)\\
&\nonumber \leq C(|\ln 2s|+\ln 3).
\end{align}
By similar arguments, the same inequality also holds true for $s\in (1,\frac{r_k}{2\rho_k})$.
Finally, we have
\begin{equation}\label{eq:subc1}
Hw_k\leq \frac{|f_k(\rho_k x,\rho_k^2t)-f_k(0)|}{\e_k}\textup{ in $S_\frac{1}{\rho_k}$},
\end{equation}
and
\begin{equation}\label{eq:cal}
Hw_k = \frac{f_k(\rho_k x,\rho_k^2t)-f_k(0)}{\e_k}\textup{ in $S_\frac{1}{\rho_k}\cap \{v_k>0\}$}
\end{equation}
with
\begin{equation}\label{eq:subcaloric}
\left(\frac{1}{r^2}\int_{S_r}\frac{(f_k(\rho_k x,\rho_k^2t)-f_k(0))^2}{\e_k^2}\d\gamma \right)^\frac12\leq \frac{1+1/k}{\mu_kC_k}\to 0,
\end{equation}
for $r\leq \frac{r_k}{\rho_k}\to \infty$, again from \eqref{eq:muk}, \eqref{eq:Msigma} and \eqref{eq:epsk}.\\

{\noindent \bf Step 2: Convergence and passing to the limit}.\\ 
Due to \eqref{eq:subquadr} for $s=2$ and \eqref{eq:subcaloric}, 
Lemma \ref{lem:compactglobal} implies that for any $\delta\in (0,1)$ there is a subsequence of $w_k$ converging in $L^2_\gamma(\R^n\times (-1,-\delta))$ 
to a limit $w_0$ and thus also in $L^2_\gammab(\R^n\times (-1,-\delta))$. By \eqref{eq:subquadr} for $s=\sqrt{\delta}\leq 1$
$$\int_{\R^n\times  \{t\geq -\delta\}}w_k^2\d\gammab\leq C\delta^2\left((\ln \delta)^2+\ln 3\right).$$
Hence, we can find a subsequence, again labelled $w_k$, converging also in $L^2_\gammab(S_1)$. By the same arguments we can also conclude that there is a subsequence converging in $L^ 2_\gammab(S_R)$ for any $R>0$. The $L^2_\gammab(S_1)$-convergence combined with \eqref{eq:T1est} and \eqref{eq:lessthanone} implies
\begin{equation}\label{eq:norm1}
\inf_{q\in \dd\mathcal{P}_{\textup{sing}}}\int_{S_1}(w_0-q)^2\d\gammab=1,
\end{equation}
with equality for $q=0$, and for all $s>1$
\begin{equation}\label{eq:normless1}
\inf_{q\in \dd\mathcal{P}_{\textup{sing}}}\int_{S_s}(w_0-q)^2\d\gammab\leq 1,
\end{equation}
where
$$\dd \mathcal{P}_{\textup{sing}}=\Big\{q=\frac12x^*Qx+mt, \tr Q=m, \exists \overline{P_k}, \frac{\overline{P_k}-P_k}{\e_k}\to q\text{ in $L^2_\gammab(S_s)$ for any $s\geq 1$} \Big\},$$
is a sort of limit differential of $\ps$. 
Moreover, passing to the limit in \eqref{eq:subc1}, \eqref{eq:subcaloric} and in \eqref{eq:cal}, using the local uniform convergence $v_k\to P_\infty$ in $\R^n\times\R^-$ (which follows from Proposition \ref{prop:quadgrowth} and Theorem \ref{thm:interior}) for some subsequential limit $P_\infty$ of $P_k$, we obtain that $w_0$ is subcaloric in $\R^n\times \R^-$ and caloric in $\R^n\times \R^-\cap \{P_\infty>0\}$. Up to a rotation of the coordinates we can assume that $\{x_1\neq 0\}\subset \{P_\infty >0\}$. Hence $w_0$ satisfies
\begin{equation}\label{eq:subcc}
Hw_0\leq 0 \text{ in $\R^n\times \R^-$},\quad  Hw_0 = 0\text{ in $\R^n\times \R^-\cap \{x_1\neq 0\}$}.
\end{equation}

\noindent {\bf Step 3: Homogeneity}.\\ 
Now we recall that Proposition \ref{prop:sing} applies for all $v_0\in \mathcal{P}_\textup{sing}$, in particular it is applicable for $v_0=P_k$ and with
 $$ 
 \bar N (r,u)=\frac{1}{r^4}\int_{S_r}(u-P_k)^2\d\gamma
 $$
so that
$$
\bar N(r,w_k)=\frac{\bar N(r\rho_k,u_k)}{\e_k^2}
$$
and
$$\frac{\d}{\d r}\bar N(r,w_k)=\frac{\rho_k}{\e_k^2} \frac{d}{\d s}\bar N(s,u_k)\big|_{s=r\rho_k}.$$
From Proposition \ref{prop:sing} and a change of variables we find that 
$$\frac{\d}{\d r}\bar N(r,w_k) \geq\frac{\rho_k}{\e_k^2}F(\rho_k r)+\frac2r\int_0^r\frac{1}{s^5}\int_{S_s}(Lw_k)^2\d\gammab \d s,$$
with $F$ as in Proposition \ref{prop:sing}, i.e.,
$$
\int_0^\tau |F(r)|\leq 4\int_0^\tau
  \frac{\sigma_2^{f_k}(r)}{r}\left(\frac{1}{r^6}\int_{B_1\times (-r^2,0]}|u_k-P_k|^2\d \gamma \right)^\frac12 \d r+32\left(\int_0^\tau
    \frac{\sigma_2^{f_k}(r)}{r} \d r\right)^2,
$$
which implies
\begin{align*}
&\int_0^R \frac{\rho_k}{\e_k^2}|F(\rho_k
r)| \d r=\int_0^{R\rho_k}\frac{|F(s)|}{\e_k^2}\d s\\ &\leq \frac{C}{\e_k^2}\left(\int_0^{R\rho_k}
  \frac{\sigma^{f_k}_2(r)}{r}\left(\frac{1}{r^6}\int_{B_1\times (-r^2,0]}|u_k-P_k|^2\d \gamma \right)^\frac12\d r+\left(\int_0^{R\rho_k}
    \frac{\sigma^{f_k}_2(r)}{r}\d r\right)^2\right).
\end{align*}
In addition, from \eqref{eq:muk}, \eqref{eq:Msigma} and \eqref{eq:epsk}, we see that $\e_k\geq \tilde C_k \Sigma(r_k)$ with $C_k = \tilde C_k (1+1/k)/\mu_k$, whenever $k$ is large enough and $\rho_k\leq \lambda_kr_k$, where
$\lambda_k\to 0$. Hence, $\e_k\geq \tilde C_k\Sigma(R\rho_k)$, for $k$ large enough. Since $\sigma^{f_k}_2\leq \sigma_p^{f_k}$ for $p\geq 2$, our choice of $\Sigma$ implies
$$
\frac{1}{\e_k}\int_0^{R\rho_k}\frac{\sigma^{f_k}_2(r)}{r}\d r\leq \frac{1}{\e_k}\int_0^{R\rho_k}\frac{\sigma^{f_k}_p(r)}{r}\d r\leq \frac{1}{\tilde C_k}\to 0.
$$
Moreover, we see that for $r\leq r_k$
\begin{align*}
\left(\frac{1}{r^6}\int_{B_1\times (-r^2,0]}|u_k-P_k|^2\d \gamma \right)^\frac12 \leq M(u_k,r_k)\leq \frac{1+1/k}{\mu_k}\e_k\leq C\e_k,
\end{align*}
due to \eqref{eq:muk} and \eqref{eq:epsk}. Hence for $R\rho_k \leq r_k$ there holds
\begin{align*}
&\frac{1}{\e_k^2}\int_0^{R\rho_k}
  \frac{\sigma^{f_k}_2(r)}{r}\left(\frac{1}{r^6}\int_{B_1\times (-r^2,0]}|u_k-P_k|^2\d \gamma \right)^\frac12\d r\\
 & \leq  \frac{C}{\e_k}\int_0^{R\rho_k}\frac{\sigma^{f_k}_2(r)}{r}\d r\leq \frac{C}{\tilde C_k}\to 0.
\end{align*}
All in all, this implies for $s_2\geq s_1$
\begin{align*}
\bar N(s_2,w_k)-\bar N(s_1,w_k)&=\int_{s_1}^{s_2} \frac{\d}{\d r}\bar N(r,w_k) \d r\\&\geq \int_{s_1}^{s_2}\frac2r\int_0^r\frac{1}{s^5}\int_{S_s}(Lw_k)^2\d\gammab \d s\d r -\int_{\rho_k s_2}^{\rho_k s_1}\frac{|F(r)|}{\e_k^2}\d r, 
\end{align*}
where the second term converges to zero as $k\to \infty$. Hence, upon passing to the limit (recall that $w_k$ converges in $L^2_\gammab(S_R)$ for any $R>0$) and using Fatou's lemma we arrive at
$$
\bar N(s_2,w_0)-\bar N(s_1,w_0)\geq \int_{s_1}^{s_2}\frac2r\int_0^r\frac{1}{s^5}\int_{S_s}(Lw_0)^2\d\gammab \d s\d r.
$$
This implies in particular for $1<R$
$$\bar N(R,w_0)-\bar N(1,w_0)\geq \int_1^R\frac2r\int_0^r\frac{1}{s^5}\int_{S_s}(Lw_0)^2\d\gammab \d s\d r.$$
Hence, from \eqref{eq:norm1} and \eqref{eq:normless1} with $q=0$ we obtain
$$0\geq \bar N(R,w_0)-\bar N(1,w_0)\geq \int_1^R\frac2r\int_0^r\frac{1}{s^5}\int_{S_s}(Lw_0)^2\d\gammab\d s\d r,$$
which forces $Lw_0=0$ in $\R^n\times (-R^2,0)$ for any $R>1$. Therefore $w_0$ is parabolically
homogeneous of degree 2. By Lemma \ref{lem:homocaloric} (recall \eqref{eq:subcc}), $w_0$ must be a caloric polynomial of degree 2.\\

\noindent {\bf Step 4: Using Liouville-type arguments towards a contradiction}.\\
By a change of coordinates we can assume that 
$$P_k=\frac12\sum_{i,j\leq n}Q^k_{ij}x_ix_j+m_kt,$$
where $\tr Q^k=m_k+1$ and $Q^k$ is a diagonal matrix, with zero entries after row
number $\ell_k$, and $m_k\in [-1,0]$.

In the set $\{t=0=x_1=\cdots=x_{\ell_k}=0\}$ we have
\begin{equation}\label{eq:nonnegk}
w_k(x,t)=\frac{u_k(\rho_k x,\rho_k^2t)}{\rho_k^2\e_k}\geq 0.
\end{equation}
Moreover, by \eqref{eq:T1est} we know that among all
$$q\in \frac{\P_\textup{sing}-P_k}{\e_k}$$
the function
$$\int_{S_1}(w_k-q)^2\d\gammab$$
attains a minimum for $q=0$. By varying $q$ we obtain
\begin{equation}\label{eq:qvariation}
\int_{S_1}w_kq\d\gammab\leq 0
\end{equation}
for all $q$ such that $\delta q\in \P_\textup{sing}-P_k$ for $\delta\geq 0$ small
enough. Notice that this set of $q$s is non-empty, due to the fact that if $q\in \ps-P_k$ then $\delta q\in \ps-P_k$ for all $\delta\in (0,1)$, since $\ps$ is convex. We can in general choose
$$q=\sum_{i=1}^n\gamma_i x_i^2+\alpha t$$
with 
$$\sum_{i=1}^n2\gamma_i=\alpha,\gamma_i\geq 0 \textup{ for $i>\ell_k$},$$
where $\alpha$ can be chosen to have any desired if $m_k\neq 0$ for all $k$ large enough (for some subsequence) (cf. {\bf Case 1} in Proposition \ref{prop:liouville})) 
and $\alpha \leq 0$ if $m_k=0$ for all $k$ large enough (cf. {\bf Case 2} in the same proposition). In this second case we
also have that for $x_1=\cdots =x_{\ell_k}=0$ and $t\leq 0$, there holds
\begin{equation}\label{eq:nonneg2}
w_k(x,t)=\frac{u_k(\rho_k x,\rho_k^2t)}{\rho_k^2\e_k}\geq 0.
\end{equation}
We also note that for $i,j\leq \ell_k, i\neq j$ and $\delta $ small enough we have
$\delta x_ix_j\in \P_\textup{sing}-P_k$. Varying $\delta $ around zero then implies
by \eqref{eq:qvariation} that
\begin{equation}\label{eq:xixjzero}
\int_{S_1}w_k x_ix_j\d\gammab=0, \textup{ for $i,j\leq \ell_k, i\neq j$}.
\end{equation}
Furthermore, with
$$R_\delta (x)=(x_1,\ldots,\underbrace{x_i\cos \delta+x_j\sin \delta}_{\textup{$i$-th coordinate}},\ldots,\underbrace{x_j\cos \delta-x_i\sin
\delta}_{\textup{$j$-th coordinate}},\ldots,x_n),$$
being the rotations in the $e_i-e_j$-plane, put
$$q_\delta(x)=P_k(R_\delta(x))-P_k(x)\in \ps-P_k,$$
for $\delta$ small enough. Therefore, by varying $\delta$ again around $\delta=0$ this implies
$$
\int_{S_1}w_k(Q^k_{ii}-Q^k_{jj})x_ix_j\d\gammab=0, \textup{for all $i,j\leq n$}.
$$
Combined with \eqref{eq:xixjzero} this implies
\begin{equation}\label{eq:rot}
\int_{S_1}w_kx_ix_j\d\gammab=0,i\neq j, i=1,\ldots,\ell_k,j=1,\ldots,n.
\end{equation}
Passing to the limit in \eqref{eq:nonnegk}, \eqref{eq:qvariation},
\eqref{eq:nonneg2} and \eqref{eq:rot} we obtain that $w_0$ satisfies
 the hypotheses of Proposition \ref{prop:liouville}. Therefore,
$w_0=0$, which is a contradiction to \eqref{eq:norm1} with $q=0$, which reads
$$
\int_{S_1}w_0^2\d\gammab = 1.
$$
\qeda\\

Now we are ready to prove Theorem \ref{thm:mainsing},
which now follows from a standard iteration, as the proof of Theorem 1.8 in \cite{Mon09}.

{\noindent \bf Proof of Theorem \ref{thm:mainsing}}\\
The proof of Theorem \ref{thm:mainsing} now follows from the
combination of Lemma 3.3, Lemma 3.4 and Lemma 3.5 in \cite{Mon09}.

\qeda\\

\section{Appendices}

\subsection{Appendix A: the connection between the different moduli of continuity}

In this section we briefly discuss the liaison between the usual $L^p$ modulus of continuity
$$\tilde{\sigma}(r)=\sup_{\rho\in(0,r]}\left(\frac{1}{|Q_\rho^-|}\int_{Q_\rho^-} |f(x,t)-f(0)|^p\d x\d t\right)^\frac{1}{p}$$
and the modulus of continuity used in this paper, namely
$$\sigma(r)=\sup_{\rho\in(0,r]}\left(\frac{1}{\rho^2}\int_{S_\rho} (f(x,t)-f(0))^p\d \gamma  \right)^\frac{1}{p}.$$
It turns out that, except for a very small error, 
$\sigma$ can be controlled by $\tilde{\sigma}(\sqrt{r})$ as proved in the proposition below.

\begin{prop}{\bf (Bound of $\sigma$ in term of $\tilde{\sigma}$)}\label{prop:modcont} \\
For any $p\in [1,+\infty)$, There are constants $c(n,p),C(n,p)>0$ such that for $r\in (0,1]$
\begin{equation}\label{eq::el5}
\sigma(r)\leq C\left(\tilde \sigma(r)+\tilde\sigma(\sqrt{ r})+\tilde\sigma(1) e^{-\frac{c}{r}}\right).
\end{equation}
In particular if $\tilde\sigma$ is Dini then so is $\sigma$.
\end{prop}

\noindent {\bf Proof of Proposition \ref{prop:modcont}}\\ 
Put $g=|f-f(0)|^p$ in order to simplify the notation during the proof. 
The proof is very similar to the one of Corollary \ref{cor:lpprim}. 
It amounts to divide the set into suitable subsets where we are able to estimate the heat kernel properly.

We begin by treating the contribution from $Q_r^-$. Let $C_k=Q^-_{2^{-k}r}\setminus Q^-_{2^{-k-1}r}$. 
From Step 2.1 of the proof of Corollary \ref{cor:lpprim}, 
for $(x,t)\in C_k$, we have 
$$|G(x,-t)|\leq \frac{C}{(2^{-k}r)^n}.$$
This gives with $r_k=2^{-k}r$
\begin{align*}
\frac{1}{r^2}\int_{Q_r^-} g\d\gamma&\leq
\frac{1}{r^2}\sum_{k=0}^\infty\int_{C_k}g\d\gamma\\
&\leq
\frac{C}{r^2}\sum_{k=0}^\infty \frac{(r_k)^{n+2}}{|Q_{r_k}^-|}\left(\int_{C_k}g\  dx dt\right) \frac{C}{(2^{-k}r)^n}\\
&\leq C\sum_{k=0}^\infty\left(\tilde\sigma(2^{-k}r)\right)^p 2^{-2k}\leq C\left(\tilde\sigma(r)\right)^p.
\end{align*}
Now we instead estimate the contribution from the part in $(B_1\setminus B_r)\times (-r^2,0]$. 
Similarly we define $\tilde{C}_k=(B_{2^{k+1}r}\setminus B_{2^kr})\times
(-r^2,0]$. 
From Step 3.1 of the proof of Corollary \ref{cor:lpprim}, there exists some $c>0$ such that for $(x,t)\in \tilde{C}_k$ we have
$$|G(x,-t)|\leq C\frac{e^{-c4^k}}{r^n},$$
Thus, we have with $r_k=2^{k+1}r$:
\begin{align*}
\frac{1}{r^2}\int_{(B_1\setminus B_r)\times
(-r^2,0]} g\d\gamma&\leq
\frac{1}{r^2}\sum_{2r \le 2^{k+1}r\le 2}\int_{\tilde{C}_k}g\d\gamma\\
&\leq \frac{C}{r^2} \sum_{2r \le 2^{k+1}r\le 2} \frac{(r_k)^{n+2}}{|Q_{r_k}^-|} \left(\int_{\tilde{C}_k}
g \ \d x\d t\right) C\frac{e^{-c4^k}}{r^n}\\
&\leq C\sum_{2r \le 2^{k+1}r\le 2}2^{k(n+2)}e^{-c4^k}\left(\tilde\sigma(2^{k+1}r)\right)^p.
\end{align*}
where we have extended $g$ by zero outside $Q_1^-$ and with the extension $\tilde{\sigma}(r)=\tilde{\sigma}(1)$ for $r\ge 1$.
Now, the last sum can be seen as a Riemann sum at points $y_k=2^kr$ so
that we obtain
\begin{align*}
\frac{1}{r^2}\int_{(B_1\setminus B_r)\times (-r^2,0]} g\d\gamma
&\leq C \sum_{k=0}^{2^{k+1}r\le 2} \frac{y_k^{n+2}}{r^{n+2}}e^{-c\frac{y_k^2}{r^2}}(\tilde\sigma(2y_k))^p\\
&\leq C \sum_{k=0}^{2^{k+1}r\le 2} \int_{y_k}^{y_{k+1}}dy\ \frac{y^{n+1}}{r^{n+2}}e^{-c\frac{y^2}{(2r)^2}}(\tilde\sigma(2y))^p\\
&\leq C\int_{r}^2dy\  \frac{y^{n+1}}{r^{n+2}}e^{-c\frac{y^2}{(2r)^2}}(\tilde\sigma(2y))^p\\
&\le C\int_{\frac12}^\frac{2}{r}ds\ s^{n+1} e^{-\frac{c}{4}s^2} (\tilde\sigma(2sr))^p
\end{align*}
This integral can be split into two parts, one when $s\leq 1/(2\sqrt{r})$
and one when $s\geq 1/(2\sqrt{r})$. For those parts we have
$$\int_{1/2}^\frac{1}{2\sqrt{r}}ds\ s^{n+1} e^{-\frac{c}{4}s^2}(\tilde\sigma(2sr))^p \leq
(\tilde\sigma(\sqrt{ r}))^p\int_{1/2}^\infty ds\ s^{n+1}e^{-\frac{c}{4}s^2}\leq C(\tilde\sigma(\sqrt{ r}))^p$$ 
and
$$\int_\frac{1}{2\sqrt{r}}^\frac{1}{r}ds\ s^{n+1}e^{-\frac{c}{4}s^2}(\tilde\sigma(2sr))^p\leq C(\tilde\sigma(1))^p \frac{1}{r^{\frac{n+2}{2}}}e^{-\frac{\bar c}{r}} \quad \mbox{with}\quad \bar c = c/16.$$
Combining the above estimate yields
$$\sigma(r)\leq C\left(\tilde\sigma(r)+\tilde\sigma(\sqrt{ r})+\tilde\sigma(1) \left(\frac{1}{r^{\frac{n+2}{2}}}e^{-\frac{\bar c}{r}}\right)^{\frac1p}\right).$$
which implies (\ref{eq::el5}).
That $\sigma$ is Dini if $\tilde\sigma$ is Dini, is immediate.
\qeda

\subsection{Appendix B: an extension theorem of Whitney-type in the parabolic setting} 

In this section, we will give a proof of a Whitney-type theorem in the
parabolic setting. For this purpose we will need a parabolic Whitney
decomposition, which is briefly described below. By a parabolic cube we here refer to a set of the form
$$Q=[x_1,x_1+r]\times \cdots\times [x_n,x_n+r]\times [x_{n+1},x_{n+1}+r^2].$$
In such a case, we define the (parabolic) ``diameter'' of the cube as
\begin{equation}\label{eq::r2}
\mbox{diam } Q = r
\end{equation}
We denote by $d$ the parabolic distance, i.e. for $(x,t)\in \R^{n+1}$ let
$$d(x,t)=\sqrt{|x|^2+|t|}$$
and we also use the notation $d(E,F)$, 
with $E,F\subset \R^{n+1}$, to denote the parabolic distance between the sets $E$ and $F$.

Given a closed set $E\subset \R^{n+1}$, what we call a Whitney decomposition is a special decomposition
of the (open) set $E^c$ into dyadic cubes. More precisely, we take a collection of cubes $\{Q_i\}$ that covers
$E^c$, with $Q_i^*$ denoting a cube with the same center as $Q_i$ but
with $\mbox{diam } Q_i^* = (1+\varepsilon)\ \mbox{diam } Q_i$ for a fixed $\varepsilon>0$.
such that
\begin{enumerate}
\item $\bigcup_{i}Q_i=E^c$,
\item the interior of the $Q_i$s are disjoint,
\item there is a constant $C$ independent of $E$ such that 
$$\frac{1}{C} \diam Q_i \leq d(Q_i,E)\leq C \diam Q_i,$$
\item there is a constant $C$ independent of $E$ such that for all $X\in Q_i^*$,
if $P_i\in E$ is a point realizing the infimum
\begin{equation}\label{eq::r3}
d(E,Q_i)=d(P_i,Q_i).
\end{equation}
Then we have
$$\frac{1}{C} d(X-P_i)\leq  \diam Q_i^*\leq C d(X-P_i),$$
and for all $A\in E$ we have
$$d(A-P_i)\leq Cd(A-X).$$
\item if $X\in E^c$ then $X$ is contained in at most a (uniform) finite number of the $Q_i^*$s.
\end{enumerate}

We also need to construct a partition of unity $\{\phi_i\}$ related to the collection $Q_i$ such that
\begin{enumerate}
\item[(i)] $\sum_i\phi_i(X)=1$ for $X\in E^c$,
\item[(ii)] $\phi_i(X)=0$ if $X\not\in Q_i^*$,
\item[(iii)] $|\nabla_{\R^n}  \phi_i(X)|\leq C(\diam Q_i)^{-1}$.
\end{enumerate}
The construction of such a decomposition together with a partition of
unity is quite standard, and can be done in the same manner as in the
case where there is no parabolic distance, see for instance p.~166-170
in \cite{Ste70}. For the parabolic case, some details can be found on
p.~188-189 in \cite{Wat02}.

To construct the collection $Q_i$ it amounts to divide the space into dyadic cubes. 
More precisely, let $M_0$ be the collection of cubes with integral corner points. 
Then we can by dividing these into smaller cubes of side length $2^{-k}$ (i.e. of diameter $2^{-k}$ with our definition (\ref{eq::r2}))
produce a sequences of collections $M_k$ that consists of the cubes with corner points in $(2^{-k}\Z^n,4^{-k}\Z)$. 
Now let $F_k$ be all the cubes in $M_k$ except those who touch $E$ or 
those who touch a cube that touches $E$. 
Let then $F=\cup_k F_k$ while ignoring all cubes that are contained in a larger cube. 
This collection of cubes will satisfy the properties above.

The construction of $\{\phi_i\}$ goes as follows. 
Take $\phi\in C_0^\infty\left(\left[-\frac{(1+\e)}{2},\frac{(1+\e)}{2}\right]^{n}\right)$ such that $\phi=1$ 
on $\left[-\frac{1}{2},\frac{1}{2}\right]^n$ and also $\psi\in C^\infty_0\left(\left[-\frac{(1+\e)}{2},\frac{(1+\e)}{2}\right]\right)$ with $\psi=1$ on $\left[-\frac{1}{2},\frac{1}{2}\right]$. Then with $X=(x,t)$,
$$\phi_i^*(X)=\phi\left(\frac{x-x^j}{\diam Q_j}\right)\psi\left(\frac{t-t^j}{(\diam Q_j)^2}\right)$$
where $(x^j,t^j)$ is the center of $Q_j$, and
$$\phi_i(X)=\frac{\phi_i^*(X)}{\sum_i\phi^*_i(X)},$$
we obtain a partition of unity satisfying our needs.

\begin{prop}\label{prop:whitney} {\bf (Parabolic $C^1$  extension of Whitney type)}\\
Take $E$ to be a closed set and assume that $f:E\subset\R^{n+1}\to \R$ and $g:E\subset\R^{n+1}\to \R^n$ satisfy
$$f(X)=f(Y)+g(Y)(x-y)+R(X,Y)$$
with $X=(x,t)$ and $Y=(y,s)$ and $|g(Y)-g(X)|\leq \omega(d(X-Y))$ and $|R(X,Y)|\leq
d(X-Y)\omega(d(X-Y))$ where $\omega$ is some monotone modulus of continuity.

Then there is a function $F:\R^{n+1}\to \R$ such that $F=f$ on $E$ which is $C^1$ regular
(for the parabolic distance $d$).
\end{prop}

\begin{rem} 
If we are given a function satisfying all the hypotheses of Proposition \ref{prop:whitney} 
except that $\omega$ is not monotone, then we can always take a slightly worse modulus of
continuity $\omega_m$ which is monotone such that the proposition applies with $\omega_m$.
\end{rem}

\noindent {\bf Proof of Proposition \ref{prop:whitney}}\\
We follow the proof of Theorem 4, p.~176 in \cite{Ste70}.
Put
\begin{equation}\label{eq::r4}
\left\{\begin{array}{lr}
F(X)=f(X)&\textup{if $X\in E$}, \\
F(X)={\sum}_{i}(f(P_i)+g(P_i)(x-p_i))\phi_i(X)& \textup{if $X\not\in E$},
\end{array}\right.
\end{equation}
where the sum is taken over all $Q_i$ and where $P_i=(p_i,t_i)\in E\subset \R^n\times \R$ is a point satisfying (\ref{eq::r3}).

\noindent {\bf Claim 1:}\\ 
For $A=(a,t_A)\in E$ we have
$$|F(X)-f(A)-g(A)(x-a)|\leq Cd(X-A)\omega(Cd(X-A)).$$
{\bf Proof of Claim 1:}\\
If $X\in E$ the claim is obvious, hence take $X\not\in E$. We compute
\begin{align*}
&|F(X)-f(A)-g(A)(x-a)|\\
&=|\sum_i(f(P_i)+g(P_i)(x-p_i)-f(A)-g(A)(x-a))\phi_i(X)|\\
&= \left|\sum_i
g(P_i)(x-p_i)-g(A)(x-a)+g(P_i)(p_i-a)-R(A,P_i))\phi_i(X)\right|\\
&=\left|\sum_i
(g(P_i)(x-a)-g(A)(x-a)-R(A,P_i))\phi_i(X)\right|\\
&\leq \sum_i
(|x-a|+d(A-P_i))\omega(d(A-P_i))\phi_i(X).
\end{align*}
Since $\phi_i$ vanishes outside $Q_i^*$ (property (ii) above), then
by property 5, the sum has only a (uniform) finite number of terms. Therefore by property 4, 
$$|F(X)-f(A)-g(A)(x-a)|\leq Cd(X-A)\omega(Cd(X-A)).$$
{\bf Claim 2:}\\
For $A\in E$ and $X\not\in E$ we have
$$|\nabla_{\R^n} F(X)-g(A)|\leq C\omega(Cd(X-A)).$$
Notice that $\nabla_{\R^n} F$ is well defined on $E^c$, because the definition (\ref{eq::r4})
only involves a finite number of terms in the sum for a given $X$ in a cube.\\
{\bf Proof of Claim 2}:\\
We remark that for $k=1,\ldots,n$ we have
$$\dd_k F(X)=\sum_i g_k(P_i)\phi_i(X)+ \sum_i (f(P_i)+g(P_i)(x-p_i))\dd_k \phi_i(X)=I_k+J_k.$$
Considering the contribution from $I_k$ we have
\begin{align*}
|I_k-g_k(A)|&\leq \sum_i|g_k(P_i)-g_k(A)||\phi_i(X)|\\
&\leq
C\sum_i\omega(d(P_i-A))|\phi_i(X)|\leq C\omega(Cd(X-A))
\end{align*}
again since the sum is finite and due to property 5 above. Furthermore,
using that property (i) implies
$$\sum\dd_k\phi_i(X)=0 \quad \mbox{for}\quad X\in E^c$$
we obtain
\begin{align*}
|J_k|&=\sum_i|f(P_i)+g(P_i)(x-p_i)-F(X)||\dd_k\phi_i(X)|\\
&\leq C\sum_i|d(P_i-X)\omega(d(P_i-X))|\frac{1}{\diam(Q_i)}\\
&\leq C\sum_i\omega(d(P_i-X))\\
&\leq C\omega(C d(A-X)).
\end{align*}
where we used property (iii) to estimate $\dd_k \phi_i$, 
Claim 1 to estimate the expansion, property (ii) and 5 to say that the sum is finite, 
and finally property 4 to relate the different distances (and $d(P_i-X)\le d(P_i-A)+d(A-X)$).
This ends the proof of claim 2.\\

The rest of the proof goes as follows: By properties (ii) and 5 the sum in
$F$ has only a finite number of terms outside $E$. Hence, $F$ is smooth
outside $E$. Moreover it has $(g,0)$ as derivative on $E$, by Claim 1. 
By the same claim it is differentiable on $E$. Therefore it is
differentiable with respect to $d$ everywhere. In addition, $g$ is continuous on $E$,
$\nabla_{\R^n}F$ is continuous on $E^c$, and by claim 2, $\nabla_{\R^n}F$ is continuous ``from $E$ to $E^c$''.
This implies that $F$ is $C^1$ in the sense of Definition \ref{def:c1}.
\qeda

\subsection{Appendix C: Justifying the partial integration}
The following result justifies the partial integration performed at several stages in the paper.

\begin{lem}\label{lem:ipp} Assume that $P$ is a homogeneous parabolic polynomial of degree $\alpha$ and that $g$ is a smooth function satisfying for $R\geq R_0>0$ and $\delta >0$
$$
\sup_{(\dd B_R)\times (-1,-\delta)}|g|\leq CR^{\beta}, 
$$
for some $\beta\geq 0$.
Then 
$$
\int_{S_1^\delta }  \nabla P \cdot \nabla g\d\gamma = -\int_{S_1^\delta }  g\nabla \cdot ((\nabla  P) G(x,-t))\d x\d t,
$$
with $S_1^\delta = \R^n\times (-1,-\delta)$.
\end{lem}
{\noindent \bf Proof of Lemma \ref{lem:ipp}} \\\noindent  Since $P$ is homogeneous of degree $\alpha$, $\nabla P$ is homogeneous of degree $\alpha-1$. Integration by parts over $B_R$ implies
\begin{align*}
&\Big|\int_{B_R\times (-1,-\delta)} ( \nabla g\cdot( \nabla P) G(x,-t)+  g \nabla \cdot ((\nabla  P) G(x,-t)))\d x\d t\Big|
\\&\leq \int_{(\dd B_R)\times (-1,-\delta)} |\nabla P||g| G(x,-t)\d\sigma \d t\\
 &\leq C\int_{(\dd B_R)\times (-1,-\delta) }R^{\alpha+\beta-1} (-t)^{-n/2}e^{\frac{R^2}{4t}}\d t\d\sigma\\
& \leq C\int_{-1}^{-\delta} R^{n+\alpha+\beta-2}(-t)^{-n/2}e^\frac{R^2}{4t} \d t.
\end{align*}
By the change of coordinates 
$$
s=\frac{R}{\sqrt{-t}},
$$
this can be bounded by
$$
C\int_R^\infty R^{\alpha+\beta}s^{n-3}e^{-s^2/4}\d s,
$$
which rapidly decays to zero as $R\to \infty$.
\qeda
\section{Acknowledgements}
Both of the authors have partially been supported by the ANR project ''MICA'', 
grant ANR-08-BLAN-0082, and by the ANR project ``HJnet'', grant ANR-12-BS01-0008-01.
Moreover, the first author has been partially supported 
by the Chair ''Mathematical modelling and numerical simulation, F-EADS - Ecole Polytechnique - INRIA'', The Royal Swedish Academy of Sciences, NTNU and MSRI. The second author has also received partial funding from the European
Research Council under the European Union's Seventh Framework Programme
(FP/2007-2013) / ERC Grant Agreement n.321186 - ReaDi -Reaction-Diffusion
Equations, Propagation and Modelling.


\bibliographystyle{amsplain}
\bibliography{ref.bib}  

\end{document}